\magnification =\magstep1    
\baselineskip =13pt    
\overfullrule =0pt    
    
\centerline {\bf ROZANSKY-WITTEN INVARIANTS VIA ATIYAH CLASSES}    
    
\vskip  .7cm    
    
\centerline{\bf M. Kapranov}
    
\vskip 1cm

Recently, L.Rozansky and E.Witten [RW] associated to any hyper-K\"ahler manifold $X$ an
invariant of topological 3-manifolds. In fact, their construction gives a system of weights    
$c_\Gamma(X)$ associated to 3-valent graphs $\Gamma$ and the corresponding    
invariant of a 3-manifold $Y$ is obtained as the sum $\sum c_\Gamma(X) I_\Gamma(Y)$    
where $I_\Gamma(Y)$ is the standard integral of the product of linking forms.    
So the new ingredient is the system of invariants $c_\Gamma(X)$ of hyper-K\"ahler    
manifolds $X$, one for each trivalent graph $\Gamma$.  They are obtained    
from the Riemannian curvature of the hyper-K\"ahler metric.     
    
\vskip .1cm

In this paper we give a reformulation of the $c_\Gamma(X)$    
  in simple cohomological terms which involve only the underlying  holomorphic symplectic    
manifold.  The idea is that we can replace  the curvature    
 by the Atiyah class [At] which is the cohomological obstruction to the existence of    
a global holomorphic  connection.    
The role of what in [RW] is called ``Bianchi identities in hyper-K\"ahler    
geometry" is played here by an identity for the square of the Atyiah class expressing    
 the existence of  the fiber bundle of second order jets. 

\vskip .1cm

The analogy between    
the curvature and the structure constants of a Lie algebra observed in [RW]    
in fact holds even without any symplectic structure, and we study the nonsymplectic case
in considerable detail so as to make the specialization to the symplectic situation easier.
We show, first of all,  that the Atiyah class of the tangent    
bundle of any complex manifold  $X$ satisfies a version of the Jacobi identity when considered    
as an element of an appropriate operad.  In particular, we find (Theorem 2.3) that for
any coherent sheaf $A$ of ${\cal O}_X$-algebras the shifted cohomology space
$H^{\bullet -1}(X, T_X\otimes A)$ has a natural structure of a graded Lie algebra,
given by composing the cup-product with the Atiyah class.  If $E$ is any holomorphic
vector bundle over $X$, then $H^{\bullet -1}(X,  E\otimes A)$ is a representation of
this Lie algebra.  

Then, we unravel the Jacobi identity to make the space of cochains
with coefficients in the tangent bundle into a  ``Lie algebra up to higher homotopies" [S].
An algebra of this type is best described by exhibiting a complex replacing the Chevalley-Eilenberg
complex for an ordinary Lie algebra. In our case this latter complex is identified with
 the sheaf of functions
on the formal neighborhood of the diagonal in $X\times X$, the identification being given by
the ``holomorphic exponential map" (the canonical coordinates of Bershadsky, Cecotti,
Ooguri and Vafa [BCOV]).  

\vskip .1cm

As far as the choice of cochains is concerned, we consider two versions. First, we use Dolbeault forms
(and assume that $X$ is equipped with a K\"ahler metric). Second, we put ourselves
into the framework of formal geometry [B] [FGG] [GGL]  and use relative forms on the space
of formal exponential maps.  The underlying algebraic result here is a 1983 theorem of
D.B. Fuks [Fuk] who described the stable cohomology of the Lie algebra of formal vector fields
with tensor coefficients in terms of what we can today  identify as the suspension of the
PROP (in the sense of  [Ad] ) governing weak Lie algebras. In the same fashion, we identify
(Theorem 3.7.4)  a certain Gilkey-type complex of natural tensors on K\"ahler manifolds
with the  suspended weak Lie PROP. 

\vskip .1cm

With the nonsymplectic case studied in detail, the introduction of a holomorphic symplectic
structure amounts to some easy modifications, presented in \S 5. 
As another outcome of our considerations    
we obtain that the $c_\Gamma(X)$ can be calculated from the curvature of an arbitrary    
K\"ahler metric, not necessarily compatible in any way with the symplectic structure.    
This may be useful because the hyper-K\"ahler metric is rarely known explicitly.     

\vskip .2cm

The outline of the paper is as follows. In \S 1 we collect some general (well known)
properties of the Atiyah classes of arbitrary holomorphic vector bundles.
In \S 2 we specialize to the case of the tangent bundle, intepret the ``cohomological
Bianchi identity" of \S 1 as the Jacobi identity and then present an unraveling of
this identity on the level of Dolbeault forms on a K\"ahler manifold.  In \S 3 we recast the properties
of the Atiyah class in the language of operads and PROPs which is well suited to treat
identities among operations such as the Jacobi identity, in an abstract way. At the end of
\S 3 we realize the weak Lie PROP by natural differential covariants on K\"ahler manifolds.
Section 4 is devoted to the formal geometry analog of K\"ahler considerations of
\S\S 2-3. Finally, in \S 5 we specialize to the case of holomorphic symplectic manifolds
and show how the previous constructions are modified and specialized in this case,
in particular, how to get the Rozansky-Witten classes $c_\Gamma(X)$ from the Atiyah
class of $X$.

\vskip .2cm    
    
The author's thinking about this question was stimulated by the letter of M. Kontsevich  [K2]     
where he sketched an interpretation of Rozansky-Witten invariants by applying    
the formalism of characteristic classes of  (symplectic) foliations to the     
the $\bar\partial$-foliation existing on $X$ considered as a $C^\infty$-manifold.     
By trying to understand his construction, the author arrived at the very elementary    
description using the Atiyah class.  However, the material of \S 4 comes closer to
Kontsevich's approach in that we use the formalism of tautological forms
familiar in the theory of characteristic classes of foliations and Gelfand-Fuks cohomology [B].

I would like to thank V. Ginzburg for showing  me Kontsevich's    
letter and  V. Ginzburg and L. Rozansky for useful discussions and suggestions
which lead the paper to evolve to its present form. Among other things, V. Ginzburg made
several suggestions about organization of the paper, in particular  that I reformulate
Theorem 2.6 in the form (2.8.1), while L. Rozansky communicated to me a formula
containing the germ of Theorem 2.6. 
I am also grateful to M. Kontsevich
who alerted me to an error in the earlier version of the text and pointed out the
reference [BCOV]. 
This work was partially supported by an NSF grant and an A.P. Sloan fellowship.

\vfill\eject    
    
\centerline {\bf \S 1. Atiyah classes in general.}    
    
\vskip 1cm

\noindent {\bf (1.1) The Atiyah class of a vector bundle.} Let $X$ be a complex analytic    
manifold (we can, if we want, work with smooth algebraic    
varieties over any field of characteristic 0). Let $E$ be a holomorphic vector bundle on $X$,    
and $J^1(E)$ be the bundle of first jects of sections of $E$. It fits into an exact    
sequence    
$$0\to \Omega^1_X\otimes E \to J^1(E)\to E\to 0,\leqno (1.1.1)$$    
which therefore gives rise to the extension class    
$$\alpha_E\in {\rm Ext}^1_X(E, \Omega^1\otimes E) = H^1(X, \Omega^1\otimes {\rm End}(E))
\leqno (1.1.2)$$    
known as the Atiyah class of $E$. An equivalent way of getting $\alpha_E$ is as follows.    
Let ${\rm Conn}(E)$ be the sheaf on $X$ whose  sections over $U\i X$ are holomorphic    
connections in $E|_U$. As well known, the space of such connections is an affine space     
over $\Gamma(U, \Omega^1\otimes {\rm End}(E))$, so ${\rm Conn}(E)$    
is a sheaf of $\Omega^1\otimes {\rm End}(E))$-torsors. Sheaves of torsors over    
any sheaf $\cal A$ of Abelian groups are classified by elements of $H^1(X,{\cal A})$,    
and $\alpha_E$ is the element classifying ${\rm Conn}(E)$. So $\alpha_E$ is an    
obstruction to the existence of a global holomorphic connection.     If $E,F$ are two vector bundles, then,
in obvious notation, we have: 
$$\alpha_{E\otimes F} = \alpha_E \otimes 1_F + 1_E \otimes \alpha_F, \leqno (1.1.3)$$
because of the well known formula for the connection in a tensor product. 

\vskip .2cm

Let ${\cal D}={\cal D}_X$ be the sheaf of rings of differential operators on $X$, and
${\cal D}^{\leq p}\i {\cal D}$ be the subsheaf of operators of order $\leq p$. It has
a natural structure of ${\cal O}_X$-bimodule, the two module structures being different.
The tensor product ${\cal D}^{\leq 1}\otimes_{\cal O} E$ is dual to $J^1(E^*)$. Therefore
$(-\alpha_E)$ is represented by the extension (symbol sequence)
$$0\to E\to {\cal D}^{\leq 1}\otimes_{\cal O} E\to T\otimes E \to 0. \leqno (1.1.4)$$

\vskip .2cm    
    
\noindent {\bf (1.2) The Bianchi identity.} If     
$a,b\in H^1(X, \Omega^1\otimes {\rm End}(E))$ are any elements,  their cup-product    
$a\cup b$ is an element of     
$H^2(X, \Omega^1\otimes \Omega^1 \otimes    
{\rm End}(E)\otimes {\rm End}(E))$. We have a natural map of vector bundles on $X$    
$$\Omega^1\otimes \Omega^1 \otimes    
{\rm End}(E)\otimes {\rm End}(E)\to S^2(\Omega^1)\otimes {\rm End}(E)    
\leqno (1.2.1)$$    
which is the symmetrization with respect to the first two arguments and the    
commutator in the second two.  We denote by $[a\cup b]\in     
H^2(X, S^2\Omega^1 \otimes {\rm End}(E))$ the image of $a\cup b$ under the map    
induced by (1.2.1) in cohomology.  

If $A,B,C$ are three sheaves on $X$ and $u\in {\rm Ext}^i(B,C)$, $v\in {\rm Ext}^j(A,B)$,
then by $u\circ v\in {\rm Ext}^{i+j}(A,C)$ we will denote their Yoneda product. 

 If $a$ is as before and $c\in H^1(X, {\rm Hom}(T\otimes T, T)) = {\rm Ext}^1(T\otimes T, T)$,
then we denote by $a*c\in H^2(X, S^2\Omega^1 \otimes {\rm End}(E))$ the Yoneda product of
the embedding $S^2T\otimes E\hookrightarrow  T\otimes T\otimes E$ and the elements
$$a\in {\rm Ext}^1(T\otimes E, E), \quad c\otimes 1\in {\rm Ext}^1(T\otimes T\otimes E, T\otimes E).$$

\proclaim (1.2.2) Proposition. The classes $\alpha_E, \alpha_{T}$ satisfy the following
property (cohomological Bianchi identity):  
$$ 2[\alpha_E\cup \alpha_E] + \alpha_E *\alpha_{T} = 0 \quad {\rm in}\quad
H^2(X, S^2\Omega^1\otimes {\rm End}(E)).$$

\noindent {\sl Proof:} Consider the two-step filtration
$$ E\i {\cal D}^{\leq 1}\otimes E\i {\cal D}^{\leq 2}\otimes E$$
with quotients $E, T\otimes E, S^2T\otimes E$  respectively.  This filtration gives the extension
classes between consecutive quotients:
$$(-\alpha_E) \in {\rm Ext}^1(T\otimes E, E), \quad \xi \in {\rm Ext}^1(S^2T\otimes E, T\otimes E),$$
whose Yoneda product is 0. Our next task is to identify $\xi$. In fact, we have the following
lemma.

\proclaim (1.2.3) Lemma. Let $s: T\otimes T\to S^2T$ be the symmetrization. Then
$\alpha_{T\otimes E} \in {\rm Ext}^1(T\otimes T\otimes E, T\otimes E)$ can be expressed
as
$$\alpha_{T\otimes E} = -\xi \circ (s\otimes 1_E) - 1_T\otimes \alpha_E.$$

The lemma implies (1.2.2) once we expand $\alpha_{T\otimes E}$ by (1.1.3).

\vskip .2cm

\noindent {\sl Proof of the lemma:} This is a particular case of a statement from [AL], n. (4.1.2.3)
applicable to any left $\cal D$-module $\cal M$ with a good filtration $({\cal M}_i)$ by
vector bundles. In such a situation we have the ``symbol multiplication" maps
$$\mu: T\otimes ({\cal M}_i/{\cal M}_{i-1})\to {\cal M}_{i+1}/{\cal M}_i$$
induced by the $\cal D$-action on $\cal M$. We also have natural extension classes
$$ f_i\in {\rm Ext}^1({\cal M}_{i+1}/{\cal M}_i, {\cal M}_i/{\cal M}_{i-1})$$

\proclaim (1.2.4) Lemma. [AL] In the described situation the class
 $\bigl(-\alpha_{{\cal M}_i/{\cal M}_{i-1}}\bigr)$ is the difference between the following
two compositions (Yoneda pairings in which the degree of {\rm Ext} is indicated by square brackets):
$$T\otimes {\cal M}_i/{\cal M}_{i-1} \buildrel \mu\over\to {\cal M}_{i+1}/{\cal M}_i
\buildrel f_i\over\to {\cal M}_{i}/{\cal M}_{i-1}[1],$$
$$T\otimes {\cal M}_i/{\cal M}_{i-1}\buildrel 1_T\otimes f_{i-1}\over\longrightarrow
T\otimes {\cal M}_{i-1}/{\cal M}_{i-2} [1]\buildrel\mu\over\to {\cal M}_i/{\cal M}_{i-1}[1].$$

To obtain Lemma 1.2.3, we take ${\cal M}={\cal D}\otimes E$ with ${\cal M}_i={\cal D}^{\leq i}\otimes E$.
Then for $i=1$ the statement identifies $(-\alpha_{T\otimes E})$. The first composition
is $\xi\circ(\sigma\otimes 1_E)$, while the second one is $-1_T\otimes\alpha_E$. This
completes the proof.

\vskip .3cm    
    
\noindent{\bf (1.3)  Atiyah class and  curvature.} The class $\alpha_E$ can be easily    
calculated both in \v Cech and Dolbeault models for cohomology. In the \v Cech    
model, we take an open covering $X=\bigcup U_i$ and pick connections    
$\nabla_i$ in $E|_{U_i}$. Then the differences    
$$\phi_{ij} = \nabla_i-\nabla_j \in\Gamma(U_i\cap U_j, \Omega^1\otimes {\rm End}(E))$$    
form a \v Cech cocycle representing $\alpha_E$.    
    
In the Dolbeault model, we pick a $C^\infty$-connection in $E$ of type $(1,0)$, i.e.,    
a differential operator    
$$\nabla: E\to \Omega^{1,0}\otimes E, \quad \nabla(f\cdot s) = \partial(f)\cdot s + f\cdot(\nabla s).$$ 
Let $\tilde\nabla=\nabla+\bar\partial$ where $\bar\partial$ is the (0,1)-connection
defining the holomorphic structure. The
 curvature $F_{\tilde\nabla}$ splits into the sum $F_{\tilde\nabla}=F_{\tilde\nabla}^{2,0}+
F_{\tilde\nabla}^{1,1}$     
according to the number of antiholomorphic differentials. Then (see [At]).    
    
\proclaim (1.3.1) Proposition.  If $\nabla$ is any smooth connection in $E$ of    
type (1,0), then $F_{\tilde\nabla}^{1,1}$ is a Dolbeault representative of $\alpha_E$.

\noindent {\bf (1.3.2) Remark.} It may be worthwhile to explain why (1.3.1) is indeed a    
{\it complete} analog of the \v Cech construction above. Namely, holomorphic    
connections in $E$ can be identified with holomorphic sections of a natural    
holomorphic fiber bundle $C(E)$, which is an affine bundle over $\Omega^1\otimes{\rm End}(E)$.    
The fiber $C(E)_x$ of $C(E)$ at $x\in X$ is the space of    
first jets of fiberwise linear isomorphisms $E_x\times X\to E$ defined near     
and identical on $E_x\times \{x\}$. Clearly, this is an affine space over $T^*_xX \otimes {\rm End}(E_x)$.    
Now,  (1,0)-connections $\nabla$ in $E$ are in natural bijection with arbitary $C^\infty$    
sections $\sigma$ of $C(E)$. Since $C(E)$ is a holomorphic affine bundle, every such    
$\sigma$ has a well defined antiholomorphic derivative $\bar\partial\sigma$ which is    
a (0,1)-form with values in the corresponding vector bundle,  i.e.,     
$$\bar\partial\sigma\in \Omega^{0,1}\otimes \Omega^{1,0}\otimes {\rm End}(E)    
=    
\Omega^{1,1}\otimes {\rm End}(E).$$    
If $\sigma$ corresponds to $\nabla$, then $\bar\partial\sigma = F_{\tilde\nabla}^{1,1}$.

\vskip .2cm

Proposition 1.3.1 has a corollary for Hermitian connections. Recall [W] that a Hermitian    
metric in a holomorphic vector bundle $E$ gives rise to a unique connection $\tilde\nabla
=\nabla+\bar\partial$    of the above type
which preserves the metric.  This connection is called the canonical connection of    
the hermitian holomorphic bundle.  It is known that $F_{\tilde\nabla}$ is in this case of type (1,1).     
Proposition 1.3.1 implies at once the following.

\proclaim (1.3.2) Proposition.     
If $E$ is equipped with a Hermitian metric and $\tilde\nabla$ is its canonical connection,     
  then $F_{\tilde\nabla}$ is a Dolbeault  representative of $\alpha_E$.     

\vskip .2cm

\noindent {\bf (1.4) Atiyah class and Chern classes.}
    If $X$ is K\"ahler, then $c_m(E)\in H^{2m}(X, {\bf C})$, the $m$th    
Chern class of $E$, can be seen as lying in $H^m(X, \Omega^m)$, and it follows    
that it is recovered from the Atiyah class as follows:    
$$c_m(E) =  {\rm Alt} ({\rm tr}(\alpha_E^{m})).\leqno (1.4.1)$$    
Here $\alpha_E^m$ is an element of $H^m(E, (\Omega^1)^{\otimes m} \otimes {\rm End}(E))$    
obtained using the tensor product in the tensor algebra  and the associative    
algebra structure in ${\rm End}(E)$, while Alt is the antisymmetrization
$(\Omega^1)^{\otimes m}\to \Omega^m$.  Note that the antisymmetrization constitutes in
fact an extra step which disregards a part of information: without it, we get an element
$$\hat c_m(E) = {\rm tr}(\alpha_E^m) \in H^m(X, (\Omega^1)^{\otimes m}). \leqno (1.4.2)$$
For a vector space $V$ let us denote by ${\rm Cyc}^m(V)$ the cyclic antisymmetric tensor power
of $V$, i.e., 
$${\rm Cyc}^m(V) = \{a\in V^{\otimes m}: ta=(-1)^{m+1} a\}, \quad t=(12...m), \leqno (1.4.3)$$
where $t$ is the cyclic permutation. Then, the cyclic invariance of the trace implies that
$$\hat c^m(E)\in H^m(X, {\rm Cyc}^m(\Omega^1)), \leqno (1.4.4)$$
but it is not, in general, totally antisymmetric. We will call $\hat c_m(E)$ the 
{\it big Chern class} of $E$; the standard Chern class is obtained from it by total antisymmetrization.

\vskip .3cm

\noindent {\bf (1.5) The Atiyah class of a principal bundle.} Let $G$ be a complex    
Lie group with Lie algebra {\bf g} and $P\to G$ be a principal $G$-bundle on $X$.     
Let ${\rm ad}(P)$ be the vector bundle on $X$ associated with the adjoint    
representation of $G$.  By considering connections in $P$, we obtain, similarly    
to the above, its Atiyah class $\alpha_P\in H^1(X, \Omega^1\otimes {\rm ad}(P))$.    
All the above properties of Atiyah classes are obviously generalized to this case.    
    
\vfill\eject    
    
\centerline {\bf \S 2.  Atiyah class of the tangent bundle and Lie brackets.}    
    
\vskip 1cm    
    
\noindent {\bf (2.1) Symmetry of the Atiyah class.}  Let $X$ be as before and $T=TX$ be    
the tangent bundle of $X$.  Specializing the  considerations of (1.1)    
to the case when $E=T$, we get the Atiyah class $\alpha_{TX}$ which    
we can  see as  an element of $H^1(X, T^*\otimes T^*\otimes T)$.     
    
\proclaim (2.1.1) Proposition.  The element $\alpha_{TX}$ is symmetric, i.e.,    
lies in the summand \hfill\break $H^1(X, S^2(T^*)\otimes T)$.    
    
\noindent {\sl Proof:}  It is enough to exhibit a sheaf of $S^2(T^*)\otimes T$-torsors    
from which  ${\rm Conn}(T)$  (a sheaf of $T^*\otimes T^*\otimes T$-torsors) is    
obtained by the change of scalars.  To find it, recall that any connection $\nabla$ in $TX$    
has a natural invariant called its {\it torsion}  $\tau_\nabla$ which is a section of     
 $\Lambda^2(T^*)\otimes T$. The sheaf ${\rm Conn}_{tf}(TX)$     
of torsion-free connections is thus a torsor    
over $S^2(T^*)\otimes T$ with required properties.    
    
\vskip .3cm    
    
\noindent {\bf (2.2) Geometric meaning of torsion-free connections.} It is convenient    
to ``materialize" the sheaf ${\rm Conn}_{tf}(TX)$ by realizing it as the sheaf of sections    
of a fiber bundle $\Phi(X)\to X$ whose fiber over $x\in X$ is an affine space over    
$S^2(T_x^*X)\otimes T_xX$. This is done as follows. For $x\in X$ let    
$\Phi_x(X)$ be the space of second jets of holomorphic maps $\phi: T_xX\to X$    
with the properties $\phi(0)=x, d_0\phi = {\rm Id}$.   A similarly defined space    
but for self-maps $T_xX\to T_xX$ is clearly just $S^2(T_x^*X)\otimes T_xX$.    
Therefore $\Phi_x(X)$ is an affine space over $S^2(T_x^*X)\otimes T_xX$.    
The  $\Phi_x(X)$ for $x\in X$ obviously unite into a fiber bundle $\Phi(X)\to X$.    
It is well known classically that sections of this bundle are the same as torsion-free    
connections.     

As a corollary of this, let us note the following interpretation of $\alpha_{TX}$
which can be also deduced from Lemma 1.2.3. 

\proclaim (2.2.1) Proposition.  The class $\alpha_{TX}$ is, up to a scalar factor,
 represented by the following
extension (second symbol sequence):
$$0\to T={\cal D}^{\leq 1}/{\cal D}^{\leq 0}\to {\cal D}^{\leq 2}/{\cal D}^{\leq 0}
\to {\cal D}^{\leq 2}/{\cal D}^{\leq 1} = S^2T\to 0.$$

We now state the first main result of this section.

\vskip .2cm

\proclaim (2.3) Theorem. Let $X$ be any complex manifold and $A$ be any quasicoherent
sheaf of commutative ${\cal O}_X$-algebras. Then:\hfill\break
(a) The maps
$$H^i(X, T\otimes A)\otimes H^j(X, T\otimes A)\to H^{i+j+1}(X, T\otimes A)$$
given by composing the cup-product with $\alpha_{TX}\in H^1(X, {\rm Hom}(S^2T,T))$,
make the graded vector space $H^{\bullet -1}(X, T\otimes A)$ into a graded Lie algebra.
\hfill\break
(b) For any holomorphic vector bundle $E$ on $X$ the maps
$$H^i(X, T\otimes A)\otimes H^j(X, E\otimes A)\to H^{i+j+1}(X, E\otimes A)$$
given by composing the cup-product with the Atiyah class $\alpha_E\in H^1(X, {\rm Hom}(T\otimes E,E))$,
make $H^{\bullet -1}(X, E\otimes A)$ into a graded $H^{\bullet -1}(X, T\otimes A)$-module.

\vskip .2cm

\noindent {\bf (2.3.1) Remarks.} (a) By construction, the structure of  a Lie algebra on the space
$H^{\bullet - 1}(X, T\otimes A)$ is bilinear over the graded commutative
algebra $H^\bullet(X,A)$, over which the former space is a module. Same for the module
structure on $H^{\bullet -1}(X, E\otimes A)$.

\vskip .1cm

(b) To see that the graded Lie algebra structure defined above is, in general, nontrivial,
it suffices to  take $A=S^\bullet(T^*)$ (the symmetric algebra with grading ignored),
$i=j=0$, and $a=b=1\in H^0(X, T\otimes T^*)$. Then the bracket $[a,b]\in H^1(X, T\otimes
S^2T^*)$ is $\alpha_{TX}$.

\vskip .1cm

(c) Theorem 2.3 is also true for sheaves of graded commutative algebras $A^\bullet$, if
we replace cohomology with the hypercohomology, i.e., consider
$${\bf H}^p(X,  T\otimes A^\bullet)=\bigoplus_{i+j=p} H^i(X, T\otimes A^j).$$

\vskip .2cm

\noindent {\bf (2.3) Proof of Theorem 2.3:} (a) If ${\bf g}^{\bullet}$ is a graded vector
space with an antisymmetric bracket $\beta: \bigwedge^2{\bf g}^\bullet \to {\bf g}^\bullet$,
then the left hand side of the Jacobi identity for $\beta$ is a  certain element
$j(\beta)\in {\rm Hom}(\bigwedge ^3 {\bf g}^\bullet, {\bf g}^\bullet)$. In our case
${\bf g}^\bullet=H^{\bullet - 1}(X, T\otimes A)$ and we find that $j(\beta)$ is given by composing
the cup product with a certain  class
$$J\in H^2(X, {\rm Hom}(S^3 T, T)).$$
This class is nothing but the symmetrization of
$$[\alpha_{TX} \cup \alpha_{TX}]\in H^2(X, S^2\Omega^1 \otimes {\rm Hom}(T,T)),$$
so it vanishes by the  ``cohomological Bianchi identity" (1.2.2) applied to $E=T$. 

\vskip .1cm

(b) If ${\bf g}^\bullet$ is a graded Lie algebra, $M^\bullet$ is a graded bector space
with a map $c: {\bf g}^\bullet\otimes M^\bullet\to M^\bullet$, then the left
hand sice of the identity
$$[g_1, g_2]m - g_1(g_2m) - (-1)^{{\rm deg}(g_1){\rm deg}(g_2)} g_2(g_1m) = 0$$
is a certain element $\tau(c)\in \bigwedge^2{\bf g}^*\otimes {\rm Hom}(M,M)$,
vanishing if and only if $M$ is a {\bf g}-module. In our case ${\bf g}^\bullet  = H^{\bullet - 1}(X, T\otimes
A)$, $M^\bullet = H^{\bullet -1}(X, E\otimes A)$, and the element $\tau(c)$ is induced
by a class
$$\sigma\in H^2(X, S^2\Omega^1\otimes {\rm Hom}(E,E))$$
which is nothing but the left hand side of (1.2.2) for $E$.  Theorem is proved.

\vskip .2cm
The case $A={\cal O}_X$ does not lead to anything interesting. Indeed, we have

\proclaim (2.3.2) Proposition. The Lie algebra structure on $H^{\bullet - 1}(X,T)$ given
by $\alpha_{TX}$, is trivial (all brackets are zero). Similarly, the module structure on
$H^{\bullet - 1}(X, E)$ is trivial.

\noindent {\sl Proof:} Let $a\in H^i(X,T)$, $b\in H^j(X,T)$. By Proposition 2.2.1, 
the bracket $[a,b]\in H^{i+j+1}(X,T)$ is obtained by applying to $a\cup b\in H^{i+j}(X, S^2T)$
the boundary homomorphism $\delta: H^{i+j}(X,S^2T)\to H^{i+j+1}(X,T)$ of the second symbol sequence.
But we have a pairing of sheaves
$$T\otimes_{\bf C}T \to {\cal D}^{\leq 1}\otimes_{\bf C} {\cal D}^{\leq 1}
\to {\cal D}^{\leq 2} \to {\cal D}^{\leq 2}/{\cal D}^{\leq 0},$$
induced by the composition of differential operators. Therefore we get an element
$a\sqcup b \in H^{i+j}(X, {\cal D}^{\leq 2}/{\cal D}^{\leq 0})$ mapping into
$a\cup b$. But this implies that $\delta(a\cup b)=0$. 

\vskip .1cm

For the bundle case the argument is similar. We note that $\alpha_E$ is represented by the
symbol sequence (1.1.4)
and that we have a pairing of sheaves
$$T\otimes_{\bf C} E\hookrightarrow {\cal D}^{\leq 1}\otimes_{\bf C} E \to 
 {\cal D}^{\leq 1}\otimes_{\cal O} E.$$

\vskip .2cm

\noindent {\bf (2.4) Weak Lie algebras and their modules.}
Theorem 2.3 is a global cohomological statement about the Atiyah class. 
We now want to give a local, cochain level strenghtening of this result.  
 Each time when the cohomology
of some sheaf forms a graded Lie algebra, it is natural to seek an underlying
structure on the space of cochains. The standard way for doing  this is by using the concept of
weak Lie algebras (or shLA's), see [S]. Let us recall this concept.

\proclaim (2.4.1) Definition. A weak Lie algebra is a {\bf Z}-graded {\bf C}-vector space
${\bf g}^\bullet$ equipped with a differential $d$ of degree $+1$ and
(graded) antisymmetric $n$-linear operations
$$b_n: ({\bf g}^\bullet)^{\otimes n} \to {\bf g}^\bullet, \quad x_1\otimes ... \otimes x_n
\mapsto [x_1, ..., x_n]_n, n\geq 2, \quad {\rm deg}(b_n) = 2-n,$$
satisfying the conditions (generalized Jacobi identities):
$$d(b_n) =\sum_{p+q=n}\sum_{\sigma\in {\rm Sh}(p,q)} {\rm sgn}(\sigma) b_{p+1}(b_q\otimes 1)
\sigma,\leqno (2.4.2)$$
where ${\rm Sh}(p,q)$ is the set of $(p,q)$-shuffles and $d(b_n)$ is the value at $b_n$ of  the natural
differential in ${\rm Hom}(({\bf g}^\bullet)^{\otimes n}, {\bf g}^\bullet)$.

In particular, $b_2$ is a morphism of complexes ($d(b_2)=0$),
 and it makes $H^\bullet_d({\bf g}^\bullet)$
into a graded Lie algebra. The higher $b_n$ are the compensating terms for the violation of
the Jacobi identity on the level of cochains rather than cohomology. 

An equivalent formulation of (2.4.2) is as follows [S]. Consider $\hat S({\bf g}^*[-1])$,
the completed symmetric algebra of the shifted dual space to {\bf g}. Each $b_n$
gives a map $b_n^*: {\bf g}^*[-1]\to S^n({\bf g}^*[-1])$ of degree 1. Let
$d_n$ be the unique odd derivation of the algebra $\hat S({\bf g}^*[-1])$ extending
$b_n^*$. Then the identities (2.4.2) all together can be expressed as one condition
$$D^2=0, \quad D=d+\sum_{n\geq 2} d_n. \leqno (2.4.3)$$
 
\vskip .1cm

Let $({\bf g}^\bullet, (b_n))$ be a weak Lie algebra. A weak ${\bf g}^\bullet$-module
is a graded vector space $M^\bullet$ equipped with a differential $d$ of degree
$+1$ and  maps
$$c_n: ({\bf g}^\bullet)^{\otimes (n-1)}\otimes M^\bullet\to M^\bullet, n\geq 2, \quad
{\rm deg}(c_n)=2-n, \leqno (2.4.4)$$
anntisymmetric in the first $n-1$ arguments and satisfying the identities which it is
convenient to express right away in the form similar to (2.4.3). Namely, let us extend
$c_n^*: M^\bullet\to M^\bullet\otimes S^{n-1}({\bf g}^*[-1])$ to a derivation $d_n^M$
of the $\hat S({\bf g}^*[-1])$-module $M^\bullet \otimes \hat S({\bf g}^*[-1])$.
Then the condition on the $c_n$ is that
$$(1\otimes D + D_M)^2=0, \quad D_M=\sum_{n\geq 2} d_n^M.\leqno (2.4.5)$$
\vskip .2cm

\noindent {\bf (2.5) Weak Lie algebra in K\"ahler geometry.}  Let $X$ be a complex manifold.
We now unravel the Jacobi identity for the Atiyah class on the level of Dolbeault forms.
Since we will work with holomorphic as well as with antiholomorphic objects, let us agree
that  in the remainder of this section
$T=TX$ will mean the holomorphic tangent bundle of $X$, while
$\Omega^{p,q}_X$ will signify the space of global $C^\infty$
forms of type $(p,q)$. Similarly,
for a holomorphic vector bundle $E$ on $X$ we will denote by $\Omega^{p,q}(E)$  the space of all $C^\infty$
forms of type $(p,q)$ with values in $E$. 

\vskip .1cm

Suppose that $X$ is equipped with a K\"ahler metric $h$.

Let $\nabla$ be the canonical (1, 0)-connection in $T$ associated with $h$, so that (1.3):
$$[\nabla, \nabla] =0 \quad {\rm in}\quad  \Omega^{2,0}( {\rm End}(T)).
\leqno (2.5.1)$$
Set $\tilde\nabla =\nabla+\bar\partial$, where $\bar\partial$ is the (0,1)-connection
defining the complex structure.  The curvature of $\tilde\nabla$ is just
$$ R  = [\bar\partial, \nabla] \in \Omega^{1,1}( {\rm End}(T)) = \Omega^{0,1}({\rm Hom}(T\otimes T,T)) \leqno (2.5.2)$$
This is a Dolbeault representative of the Atiyah class $\alpha_{TX}$, in particular,
$$\bar\partial R = 0 \quad {\rm in}\quad \Omega^{0,2}({\rm Hom}(T\otimes T, T)) \leqno (2.5.3)$$
(Bianchi identity). Further, the condition for $h$ to be K\"ahler is equivalent, as it is well known,
to torsion-freeness of $\nabla$, so actually
$$R\in \Omega^{0,1}( {\rm Hom}(S^2T,T)). \leqno (2.5.4)$$
Let us now define tensor fields $R_n$, $n\geq 2$, as higher covariant
derivatives of the curvature:
$$R_n\in \Omega^{0,1}({\rm Hom}(S^2T\otimes T^{\otimes(n-2)}, T)), \quad
R_2:=R, \, R_{i+1}=\nabla R_i.\leqno (2.5.5)$$

\proclaim (2.5.6) Proposition. Each $R_n$ is totally symmetric, i.e.,
$R_n\in \Omega^{0,1} ({\rm Hom}(S^nT,T))$.

\noindent {\sl Proof:} Follows immediately from (2.5.1).

\vskip .2cm

 Except for $R_2=R$ the forms $R_n$ are not, in general, $\bar\partial$-closed.
Let $\Omega^{0,\bullet}(T)$ be the Dolbeault complex of global smooth $(0,i)$-forms
with values in $T$, and $\Omega^{0, \bullet-1}(T)$ be the shifted complex.

\vskip .1cm

\proclaim (2.6) Theorem. The maps
$$b_n: \Omega^{0, j_1}(T)\otimes ... \otimes \Omega^{0, j_n}(T)
\to \Omega^{0, j_1+...+j_n+1}(T), \quad n\geq 2,$$
given by composing the wedge product (with values in $\Omega^{0,\bullet}(T^{\otimes n})$)
with \hfill\break
$R_n\in \Omega^{0,1}({\rm Hom}(T^{\otimes n}, T))$, make the
shifted Dolbeault complex $\Omega^{0, \bullet-1}(T)$ into a weak Lie algebra. 

\proclaim (2.6.1) Corollary. If $X$ is a Hermitian symmetric space, then $R$ makes
$\Omega^{0, \bullet -1}(T)$ into a genuine Lie dg-algebra.

\noindent{\sl Proof:} We need to establish the generalized Jacobi identities (2.4.1) for the $R_n$.
For this, write:
$$\bar\partial R_n = \bar\partial \nabla ... \nabla R \leqno (2.6.2)$$
(with $(n-2)$ instances of $\nabla$) and use the  commutation relation (2.5.2)
 together with (2.5.3).
 This gives
$$\bar\partial R_n= \sum_{a+b=n-2} \nabla^a \circ R_*\circ \nabla^bR, \leqno (2.6.3)$$
where
$$R_*\in \Omega^{1,1} ({\rm End}({\rm Hom}(S^{b+2}T,T)))$$
is the operator-valued (1,1)-form induced by $R$.  By evaluating $R_*$, we find
$$\bar\partial R = \sum_{p+q=n}\sum_{\sigma\in {\rm Sh}(p,q)} R_{p+1}(R_q\otimes 1)\sigma,$$
which  differs from the RHS of the generalized Jacobi identity only by the
absense of the signs ${\rm sgn}(\sigma)$. These
signs, however, constitute exactly the effect of shift from $H^\bullet$ to $H^{\bullet-1}$. Theorem is 
proved. 

\vskip .2cm

\noindent {\bf Remark.} The first instance of Theorem 2.6 (that $R_3$ cobounds the Jacobi identity for the curvature) was communicated to me by L. Rozansky. 

\vskip .3cm

\noindent {\bf (2.7) Companion theorem for vector bundles.} Let now $(E, h_E)$ be a Hermitian 
holomorphic vector bundle on a K\"ahler manifold $X$, and let $\nabla_E$ be its canonical
(0,1)-connection, so that
$$ [\nabla_E, \nabla_E]=0 \quad {\rm in}\quad \Omega^{2,0}({\rm End}(E)).\leqno (2.7.1)$$
Let
$$F= [\bar\partial, \nabla_E]\in \Omega^{1,1}({\rm End}(E)) = \Omega^{0,1}({\rm Hom}(T\otimes
E,E))\leqno (2.7.2)$$
be the total curvature of $\nabla_E$. Then
$$\bar\partial F=0 \quad {\rm in}\quad \Omega^{2,0}({\rm Hom}(T\otimes E,E)), \leqno (2.7.3)$$
and $F$ is the Dolbeault representative of the Atiyah class $\alpha_E$. Define the tensor fields
$$F_n\in\Omega^{0,1}({\rm Hom}(S^{n-1}T\otimes E,E))\leqno (2.7.4)$$
by setting
$$F_2=F, \,\, F_n=\nabla F_{n-1}, n\geq 3.\leqno (2.7.5)$$
As before, the required symmetry of $F$ follows from (2.7.1). 

\proclaim (2.7.6) Theorem. The maps
$$c_n: (\Omega^{0, \bullet -1}(T))^{\otimes (n-1)}\otimes \Omega^{0, \bullet -1}(E)\to
 \Omega^{0, \bullet -1}(E)$$
given by composing the wedge product with $F_n$, make the Dolbeault complex
$ \Omega^{0, \bullet -1}(E)$ into a weak module over the weak Lie algebra
$ \Omega^{0, \bullet -1}(T)$.

The proof, using (2.7.1-3), is almost identical to that of (2.6) and is left to the reader. 

\proclaim (2.7.7) Corollary. If $(E, h_E)$ is a homogeneous Hermitian bundle over a
Hermitian symmetric space $X$, then $F$ makes $\Omega^{0, \bullet - 1}(E)$ into
a dg-module over the dg-Lie algebra $\Omega^{0, \bullet - 1}(T)$.

\vskip .2cm

\noindent {\bf (2.8) Interpretation via $D^2=0$.} In the notation of  (2.5),
let
$$R_n^* \in\Omega^{0,1}({\rm Hom}(T^*, S^nT^*))$$
be the partial transpose of $R_n$. Consider the completed symmetric algebra $\hat S(T^*)$
(this is a sheaf of ungraded ${\cal O}_X$-algebras) and introduce in the algebra
$\Omega^{0,\bullet}(\hat S(T^*))$ the grading induced from that on $\Omega^{0,\bullet}$.
Let $\tilde R_n^*$ be the odd derivation of this algebra induced by $R_n^*$.
Theorem 2.6 can be reformulated as follows (I am grateful to V. Ginzburg for suggesting that
I do this). 

\proclaim (2.8.1) Reformulation.  The derivation   $D=\bar\partial +\sum_{n\geq 2} \tilde R_n^*$ of
$\Omega^{0,\bullet}(\hat S(T^*))$ satisfies $D^2=0$. 

This is not exactly the result of applying (2.4.3) to ${\bf g}^\bullet = \Omega^{0, \bullet -1}(T)$ because we
take symmetric powers over ${\cal O}_X$ rather than {\bf C} and also do not seem to dualize
the spaces $\Omega^{0,\bullet}$.  But because the maps $R_n$ are ${\cal O}_X$-linear
and because the formal adjoint of $\bar\partial: \Omega^{0,i}\to \Omega^{0, i+1}$
is $\bar\partial: \Omega^{0, r-i-1}\to \Omega^{0,r-i}$ ($r={\rm dim}(X)$), this
change of context is justified. 

\vskip .2cm

Let us now view $D$ geometrically. The sheaf $\hat S(T^*)$ is the sheaf of functions on
$X^{(\infty)}_{TX}$, the formal neighborhood of $X$ (regarded as the zero section)
in (the total space of ) $TX$. More formally, denoting by $\pi: X^{(\infty)}_{TX}\to X$ the natural
projection, we can write
$$\hat S(T^*) = \pi_*\bigl({\cal O}_{X^{(\infty)}_{TX}}\bigr).$$
The derivation $D$ in $\Omega^{0,\bullet}(\hat S(T^*))$ can thus be regarded as a non-linear
(0,1)-connection {\bf D}  in the fiber bundle $\pi: X^{(\infty)}_{TX}\to X$. The condition
$D^2=0$ means that {\bf D} is integrable, i.e., defines a new holomorphic structure in
$X^{(\infty)}_{TX}$.  We are going to describe this new structure and at the same time  give
a very natural explanation of the previous constructions. Namely, consider $X^{(\infty)}_{X\times X}$,
the formal neighborhood of the diagonal $X\i X\times X$. This is a fiber bundle over $X$
(with respect to the projection to, say, the second factor) whose fiber over $x\in X$
is $x^{(\infty)}_X = {\rm Spf}(\hat {\cal O}_{X,x})$, the formal neighborhood of $x$ in $X$. Clearly, 
this fiber bundle has a holomorphic structure induced from that on $X$.

\proclaim (2.8.2) Theorem. The bundle $X^{(\infty)}_{TX}$ with the new complex structure {\bf D}
is naturally isomorphic to $X^{(\infty)}_{X\times X}$.

The proof is given in the next subsection.

\vskip .3cm

\noindent {\bf (2.9) The holomorphic exponential map.} Let $x\in X$ be a point.  Recall that by
$T_xX$ we denote $T_{x}^{1,0}X$, the ``holomorphic" tangent space which we want to
distinguish from $T_x^{\bf R}X$, the tangent space to $X$ conisdered as a real manifold.
More precisely, let $I: T_x^{\bf R}X\to T_x^{\bf R}X$ be the complex structure, $I^2=-1$,
and $T_x^{\bf C}X = {\bf C}\otimes_{\bf R}T_x^{\bf R}X$. Then $T_xX$ is the $(+i)$-eigenspace
of $1\otimes I$ on  $T_x^{\bf C}X$. The correspondence $\xi\mapsto \xi-iI\xi$ defines
an isomorphism of complex vector spaces $(T_x^{\bf R}X, I)\to
(T_xX, i)$. 

Now, the geodesic exponential map at $x$ (for $X$ considered as a real manifold)
$$\exp_x^{\bf R}: T_x^{\bf R}X \to X$$
is not, in general, holomorphic. Suppose first that our K\"ahler metric is real analytic.
Then so is $\exp_x^{\bf R}$, and we can take its analytic continuation ``to the complex
domain". 
 In other words, let $X'=X$ and $X''$ be $X$
with the opposite complex structure. Then the image of the diagonal embedding
$X\hookrightarrow X'\times X''$ is totally real, so $X'\times X''$ can be seen as the
complexification of $X$. Therefore $\exp_x^{\bf R}$ continues to a holomorphic map
$$\exp_x^{\bf C}: T_x^{\bf C}X = T_xX \oplus T_x^{0,1}X \to X'\times X'',$$
defined in some neighborhood of 0.

\proclaim (2.9.1) Lemma. 
Suppose the K\"ahler metric on $X$ is real analytic. Then,  the restriction of $\exp_x^{\bf C}$ to $T_xX$ takes values in
$X'\times \{x\}$ and thus gives (via the holomorphic identification $X'\to X$) a holomorphic map
$\exp_x: T_xX\to X$ defined in some neighborhood of 0, and whose differential at 0 is the
identity.

\noindent {\sl Proof:}  The complexified Riemannian connection on $T^{\bf C}X$ is $\tilde\nabla =
\nabla +\bar\partial$. Its analytic continuation is a holomorphic connection 
$\tilde\nabla^{\bf C} = \nabla^{\bf C} = \bar\partial^{\bf C}$ in the holomorphic
tangent bundle of $X'\times X''$, defined in some neighborhood of $X$.  The summands $\nabla^{\bf C}$
and $\bar\partial^{\bf C}$ have types (1,0) and (0,1) with respect to the
decomposition
$$T_{(x', x'')}(X'\times X'') = T_{x'}X' \oplus T_{x''} X''.$$
This decomposition being flat for $\nabla^{\bf C}$ and holomorphic, the exponential map for $\nabla^{\bf C}$
at a diagonal point $(x,x), x\in X$ takes $T_xX'$ into $X'\times \{x\} \simeq X$. But $T_xX'
\i T_{(x,x)}(X'\times X'')$ is precisely $T_x^{1,0}X\i T_x^{\bf C}X$, and the exponential map
for $\nabla^{\bf C}$ is just the restriction of $\exp_x^{\bf C}$  to $T_x^{1,0}$. Lemma is proved.

\vskip .2cm

The map $\exp_x$ can be called the {\it holomorphic exponential map}. It seems that it was
first introduced in 1994, in the physical paper [BCOV], \S 2.6 and called ``canonical coordinates".
Note that even when the metric is not analytic but only smooth, consideration of the Taylor expansion
of $\exp_x^{\bf R}$ in coordinates $z_i, \bar z_i$ (where $z_i$ form a local holomorphic
coordinate system), furnishes an isomorphism of formal neighborhoods
$$\exp_x: 0^{(\infty)}_{T_xX}\to x^{(\infty)}_{X}, \leqno (2.9.2)$$
which will be sufficient for the purposes we have in mind.

\vskip .2cm

\noindent {\bf (2.9.3) Example. }  Let $X= {\bf C}P^1$ with the Fubini-Study metric. As a Riemannian
manifold, $X$ is the unit sphere $S^2\i {\bf R}^3$.  Choose a point $x\in X$ and introduce
in $T_xX$ a linear coordinate system $(u,v)$ be means of an orthogonal frame. Then identify
a neighborhood of $x$ with $T_xX$ by means of the stereographic projection from the opposite
point, thereby introducing a coordinate system in $X$. Elementary trigonometry gives:
$$\exp_x^{\bf R}(u,v) = {2\sin \sqrt{u^2+v^2}\over
\sqrt{u^2+v^2} (1+\cos \sqrt{u^2+v^2})} (u,v)$$
(this is real analytic since $\sin(z)/z$ and $\cos(z)$ are even functions).  Now, the complex structure
in $T_x^{\bf R}X$ is $I(a,b) = (b, -a)$. 
Thinking now of $u,v$ as complex
variables and substituting $u=a-ib, v=b+ia$ with $a,b\in {\bf R}$ (which means that we restrict to
$T_x^{1,0}\i T_x^{\bf C}$) we find that the radicals vanish and we get
$\exp_x(z) = z$, $z\in T_xX$. So the holomorphic exponential map is, in this case,
exactly the stereographic projection, i.e., the affine coordinate on ${\bf C}P^1$ for which
the point opposite to $x$ serves as the infinity.  In a similar way, for  $X$ a Grassmannian the
map $\exp_x$ provides an affine identification of $T_xX$ with an open Schubert cell. 

\vskip .2cm

Let us now prove Theorem 2.8.2. Consider, for any $x\in X$, the formal isomorphism (2.9.2). 
These isomorphisms unite into a fiberwise holomorphic isomorphism of fiber bundles
$$\exp: X^{(\infty)}_{TX}\to X^{(\infty)}_{X\times X}.$$
The variation with respect to $x$ of the $\exp_x$ is not, in general, holomorphic in the
usual sense. However, we have the following statement which implies our theorem.

\proclaim (2.9.4) Proposition.  The map $\exp$ is holomorphic with respect to the complex
structure {\bf D} on $X^{(\infty)}_T$.

\noindent {\sl Proof:}  We will consider the real analytic case. The general case presents
only notational complications in that we replace $X'$ and $X''$ below by working in the
variables $z_i$ and $\bar z_i$.

By considering the connection $\nabla^{\bf C}$ on $X'\times X''$,
we reduce ourselves to the following purely holomorphic problem.

\vskip .2cm

Suppose given a complex manifold $X$ and a family $\nabla = (\nabla_s)_{s\in S}$
of flat torsion-free connections in $TX$ parametrized by some complex manifold $S$.
Let $p,q$ be the projections of $X\times S$ to $X$ and $S$ respectively. Then the variation
(derivative) of the $\nabla_s$ with respect to $s$ is a section
$$R\in \Gamma(X\times S, \, q^*\Omega^1_S \otimes p^*{\rm Hom}(S^2TX, TX)).$$
We can apply to each restriction $R|_{X\times \{s\}}$ the covariant derivative $\nabla_s$ several times,
getting tensor fields
$$R_n = \nabla^{n-2}R \in \Gamma(X\times S, \, q^*\Omega^1_S \otimes
p^* {\rm Hom}(S^nTX, TX)). \quad n\geq 2.$$
On the other hand, for every $x,s$ the connection $\nabla_s$ gives rise to the exponential
map
$$\exp_{x,s}: T_xX\to X, \quad 0\mapsto x,  d_0 \exp_{x,s} = {\rm Id},$$
whose variation with respect to $s$ is, for each fixed $x$, a 1-form on $X$ with values in analytic vector
fields on (some neighborhood of 0 in) $T_xX$ with vanishing constant and linear terms.
Recall that for any vector space $V$ the space of formal vector fields on $V$ at 0 is the
product $\prod_{n\geq 0} {\rm Hom}(S^nV, V)$. Thus we can write the Taylor expansion of the
variation as
$$\exp_{x,s}^{-1} d_s \, \exp_{x,s} \in \Gamma \biggl(X\times S, \, q^*\Omega^1_S\otimes \prod_{n\geq 2}
p^*{\rm Hom}(S^nTX, TX)\biggr).$$
In order to establish Proposition 2.9.4, it is enough to prove the following.

\proclaim (2.9.5) Proposition. In the described situation $R_n$ is the $n$th homogeneous component
of $\exp_{x,s}^{-1} d_s \, \exp_{x,s}$.

\noindent {\sl Proof:} Fix some $x_0\in X, s_0\in S$ and identify $T_{x_0}X$ with ${\bf C}^r$,
$r=\dim(X)$ by means of some linear isomorphism. Then use $\exp_{x_0, s_0}$
as a coordinate system on $X$ near $x_0$. For any $s$ the connection $\nabla_s$
is then defined in our coordinates by its connection matrix $\Gamma(s)
\in \Gamma({\bf C}^r, {\rm Hom}(S^2T, T))$, so that $R=d_s\Gamma(s)$ is just its derivative
with respect to $s$. For $s=s_0$ we have $\Gamma(s_0)=0$, because the exponential map for
a flat torsion free connection takes it into the  standard Euclidean connection on the tangent space.
This implies that the higher covariant derivatives $\nabla_{s_0}^i R|_{X\times\{s_0\}}$
are the same as the usual derivatives, with respect to our chosen coordinates, of
$R_{s_0} = d_s|_{s=s_0} \Gamma(s)$.  By the same token as before,  for arbitrary $s$ the flatness of
$\nabla_s$ allows us to describe it as the connection induced from the standard Euclidean
connection on ${\bf C}^r$ by the change of coordinates given by $\exp_{x_0, s}$. 
So our statement  reduces to the following lemma.

\proclaim (2.9.6) Lemma. Let $v = \sum_{i=1}^r v_i(z) \partial/\partial z_i$ be
a holomorphic vector field on (some domain of) ${\bf C}^r$. Regarding $v$ as an infinitesimal
diffeomorphism (i.e., the tangent to a family of diffeomorphisms $g(s): {\bf C}^r\to {\bf C}^r$,
$s\in {\bf C}$, $g(0)={\rm Id}$), let $\Gamma\in\Gamma({\bf C}^r, {\rm Hom}(S^2T,T))$
be the correspodning infinitesimal  variation of the connmections (induced by the $g(s)$
from the Euclidean one).  Then the components of $\Gamma$ are
$$\Gamma^i_{jk}(z) = {\partial^2 v_i\over \partial z_j\partial z_k}.$$

The proof of this lemma is straightforward from the standard formulas of differential
geometry.

\vfill\eject

\centerline {\bf \S 3. Operadic interpretation.}

\vskip 1cm

As we saw, for any sheaf $A$ of commutative algebras on $X$,
the Atiyah class $\alpha_{TX}\in H^1(X, {\rm Hom}(S^2T, T))$ makes
each $H^{\bullet-1}(X, T\otimes A)$, into a graded Lie algebra. Each composite $m$-ary operation
in this algebra (such as, e.g., $[[x_1, x_2], [x_3, x_4]]$ for $m=4$) is represented by
a certain class in $H^{m-1}(X, {\rm Hom}(T^{\otimes m}, T))$
composed out of $\alpha_{TX}$. In this section we study these classes by themselves rather
than by using the operations on $H^{\bullet -1}(X, T\otimes A))$ represented by them.
For this, we use the language of operads and PROPs, see   [Ad] [GiK], [GeK1-2] [KM].

\vskip .3cm

\noindent {\bf (3.1) Reminder on operads, PROPs and modules.}
Recall that an operad $\cal P$ is a collection of vector spaces ${\cal P}(n)$, $n\geq 0$,    
together with the action of $S_n$, the symmetric group, on ${\cal P}(n)$ for each $n$    
and composition maps     
$$\circ_i: {\cal P}(m)\otimes {\cal P}(n)\to {\cal P}(m+n-1), \quad i=1, ..., m.$$    
satisfying appropriate equivariance and associativity axioms.  Informally,    
elements of ${\cal P}(m)$ can be thought as $m$-ary operations, the $S_m$-action    
as permutation of arguments in the operations, and $p \circ_i q$ as the     
operation     
$$p(x_1, ..., x_{i-1}, q(x_i, ..., x_{i+n-1}), x_{i+n}, ..., x_{m+n-1}) \leqno (3.1.1)$$    
    
An algebra over an operad    
$\cal P$ is a vector space $A$ together with $S_n$-invariant maps $\mu_n: {\cal P}(n)\otimes    
A^{\otimes n}\to A$ satisfying the associativity properties which mean that the    
compositions $\circ_i$ in $\cal P$  indeed go, under the $\mu_n$, into the substitution    
of one operation inside another, as described in (3.1.1).

\vskip .1cm

The concept of a PROP (see [Ad]) is slightly more general.      
While operads
describe algebras $A$ with operations of the form $A^{\otimes n}\to A$, PROPs allow for
more general operations $A^{\otimes n}\to A^{\otimes m}$ (which may or may not be
deducible from the former ones). 

Thus a PROP  $\Pi$ is a family of vector spaces $\Pi(n,m)$, $n,m\geq 0$, equipped
with a left $S_n$-action and a right $S_m$-action, commuting with each other, as well
as the following structures:

\vskip .1cm

\item{(3.1.2)} Composition maps $\Pi(n,p)\otimes \Pi(m,n)\to \Pi(m,p)$, making $\Pi$ into 
a category with the set of objects $[m], m\in {\bf Z}_+$ and ${\rm Hom}([n], [m]) =
\Pi(n,m)$.

\vskip .1cm

\item{(3.1.3)} Juxtaposition maps $\Pi(n,m)\otimes \Pi(n',m')\to \Pi(n+n', m+m')$, making
$\Pi$ into a symmetric monoidal category with monoidal operation on objects defined
by $[n]\odot [n'] = [n+n']$.

\vskip .2cm

A (right) module over an operad $\cal P$ is (see [M]) a collection $\cal M$ of $S_n$-modules
${\cal M}(n)$, $n\geq 0$ and compositions
$$\circ_i: {\cal M}(m)\otimes {\cal P}(n)\to {\cal M}(m+n-1), i=1, ..., m$$
satisfying the equivariance and associativity axioms obtained by polarizing those of an operad.

\vskip .2cm

\noindent {\bf (3.1.4) Examples.} (a)
For any vector space $V$ we have its {\it endomorphism operad} ${\cal E}_V$    
with components ${\cal E}_V(n) = {\rm Hom}(V^{\otimes n}, V) = (V^*)^{\otimes n}\otimes V$.    
The space $V$ is canonically an algebra over ${\cal E}_V$.  For any operad $\cal P$    
a structure of $\cal P$-algebra on a vector space $A$ is the same as a morphism of    
operads ${\cal P}\to {\cal E}_A$.   

Similarly,    we have a PROP 
${\rm END}_V$ with ${\rm END}_V(n,m) = {\rm Hom}(V^{\otimes n}, V^{\otimes m})$.
An algebra over a PROP  $\Pi$ is a vector space $A$ together with a morphism of PROPs
$\Pi\to {\rm END}_A$. For example,  the class of Hopf algebras can be described by a PROP
 but not an operad.

\vskip .1cm

(b) Any operad is a module over itself. If $\Pi$ is a PROP, then the spaces ${\cal P}(n) = \Pi(n,1)$
form an operad. For every $k$ the spaces $\Pi_a(n) = \Pi(n,a)$
form a module over this operad. 

\vskip .3cm

\noindent {\bf (3.2) dg-operads and PROPs. } All the above constructions can be carried out in any symmetric
monoidal category. 
By a differential graded (dg-) operad we mean an operad in the symmetric monoidal category of 
differential graded vector spaces, i.e., cochain complexes
(in that category the symmetry isomorphisms are given by the Koszul sign rule). For a  dg-vector space    
$V^\bullet$ we define its shifts $V^\bullet[m]$ by $(V^\bullet[m])^i = V^{m+i}$.     
For a dg-operad $\cal P$ its {\it suspension} $\Sigma ({\cal P})$    
is a new dg-operad formed by the shifted spaces $\Sigma ({\cal P}) (n) = {\cal P}(n)[1-n]$    
with the symmetric group action differing from that on ${\cal P}(n)$    
be tensoring with the sign representation, see [GeK1] for the explicit formulas for    
the compositions.     
 If $A^\bullet$ is a differential graded $\cal P$-algebra,    
then  ${\cal A}[1]$ is a $\Sigma ({\cal P})$-algebra. For $p\in {\cal P}(n)$    
let $\Sigma (p)$ be the corresponding element of $\Sigma ({\cal P})(n)$.   
The conventions for PROPs are similar. Thus, the suspension $\Sigma\Pi$ of a dg-PROP $\Pi$
has $\Sigma\Pi(n,m) = \Pi(n,m)[m-n]$. We will view graded vector spaces as dg-vector
spaces with zero differential. 

\vskip .3cm

\noindent {\bf (3.3) A PROP from an operad.} 
Let $\cal P$ be an operad. 
We  define a $\cal P$-module ${\cal P}(-,0) = \{ {\cal P}(n,0)\}$ called the
module of natural forms (on $\cal P$-algebras). It is defined as the $\cal P$-module
generated by symbols
$${\rm tr}(p)\in {\cal P}(n,0), \quad p\in {\cal P}(n+1), \leqno (3.3.1)$$
subject to the following relations:
$${\rm tr}(p\sigma) = {\rm tr}(p)\sigma, \quad \sigma\in S_n\i S_{n+1}, \leqno (3.3.2)$$
$${\rm tr}(p \circ_i q) = {\rm tr}(p) \circ_i q, \quad p\in {\cal P}(a+1),
q\in {\cal P}(b+1), i\neq a+1, \leqno (3.3.3)$$
$${\rm tr}(p \circ_{a+1}q) = tr(q \circ_{b+1}p)\tau, \quad
\tau=\pmatrix{1&2&...&a&...a+b\cr a+1&...&a+b&1&...&a}, \leqno (3.3.4)$$
$p\in {\cal P}(a+1), q\in {\cal P}(b+1).$

\vskip .2cm

Motivation: if $A$ is a finite-dimensional $\cal P$-algebra, then any $p\in {\cal P}(n+1)$
gives a morphism $\mu_p: A^{\otimes (n+1)}\to A$, and we can take its trace
${\rm tr}_{n+1}(\mu_p): A^{\otimes n}\to {\bf C}$ with respect to the last contravariant
argument and the only covariant argument. The requirements on the ${\rm tr}(p)$ are the
axiomatizations of the properties of these traces.

\vskip .1cm

We now define a PROP, denoted $\Pi_{\cal P}$ to be generated by formal juxtapositions
and permutations from ${\cal P}(n,0)=\Pi_{\cal P}(n,0)$ and
${\cal P}(n)\i {\Pi}_{\cal P}(n,1)$. In other words,
$$\Pi_{\cal P}(n,m) = \bigoplus_{\{1, ..., n\}
=\atop A_1\cup ... \cup A_m\cup B_1\cup...\cup B_r}
\bigotimes_i {\cal P}(A_i) \otimes\bigotimes_j {\cal P}(B_j,0), \leqno (3.3.5)$$
where ${\cal P}(A), \#(A)=a$, is the notation for the functor on the category of
$a$-element sets and their bijections associated to the $S_a$-module ${\cal P}(a)$. 

\proclaim (3.3.6) Proposition. If $A$ is a finite-dimensional ${\cal P}$-algebra, then it is
also a $\Pi_{\cal P}$-algebra.

\vskip .2cm

\noindent {\bf (3.4) The Lie operad and PROP.} 
 We denote by     
${\cal L}ie$ the Lie operad, whose algebras are Lie algebras in the usual sense,     
see [GeK1-2] [GiK].  Explicitly, ${\cal L}ie(n)$ is a subspace in the free Lie algebra on generators    
$x_1, ..., x_n$ spanned by Lie monomials containing each $x_i$ exactly once.    
Thus ${\cal L}ie(2)$ is one-dimensional and spanned by $[x_1, x_2]$    
(which is anti-invariant under $S_2$), while ${\cal L}ie(3)$ is two-dimensional and spanned by    
three elements    
$$[x_1, [x_2,x_3]], \quad [x_2, [x_1, x_3]], \quad [x_3, [x_1, x_2]]$$    
whose sum is zero (Jacobi identity).    
Given an arbitrary operad $\cal P$    
and an element $p\in {\cal P}(2)$, we will say that $p$ is a {\it Lie element}, if $p$    
is antisymmetric and satisfies the Jacobi identity. In other words, $p$ is a Lie element    
if there is a morphism of operads ${\cal L}ie\to {\cal P}$ which  takes the generator    
$[x,y]\in {\cal L}ie(2)$ into $p$. Such a morphism is unique, if it exists.    

\vskip .2cm

 We denote the PROP $\Pi_{{\cal L}ie}$ by LIE. The new
generators in LIE (apart from the bracket $[x_1, x_2]\in {\rm LIE}(2,1)$) form the space
${\rm LIE}(n,0) = {\cal L}ie(n,0)$.  An example of an element of the latter space is
given by
$$\kappa_n = {\rm tr}([x_1 [x_2 ...[x_n, x_{n+1}]...]).\leqno (3.4.1)$$
Here $[x_1 [x_2 ...[x_n, x_{n+1}]...]$ is regarded as an element of ${\cal L}ie(n+1)$.
It follows from (3.3.4) that $\kappa_n$ is cyclically symmetric, i.e.,
$$\kappa_n t = \kappa_n, \quad t=(12...n)\in {\bf Z}_n\i S_n.\leqno (3.4.2)$$
If {\bf g} is a finite-dimensional Lie algebra, then $\kappa_n$ gives the $n$th Killing form
on {\bf g}:
$$x_1\otimes ... \otimes x_n \mapsto {\rm tr}( {\rm ad}(x_1) ... {\rm ad}(x_n)).
\leqno (3.4.3)$$

\proclaim (3.4.4) Proposition. The space ${\cal L}ie(n,0)$ has dimension $(n-1)!$ and a basis there
is formed by the elements $\kappa_n\sigma$, $\sigma\in S_n/{\bf Z}_n$. 

\noindent {\sl Proof:} This follows from the fact that a basis in ${\cal L}ie(n+1)$
is formed by the Lie  monomials
$$[x_{\sigma_1}[ ...[x_{\sigma(n)}, x_{n+1}]...], \quad \sigma\in S_n.$$

\vskip .3cm

\noindent {\bf (3.5) The Atiyah class as a Lie element.} Let Now $X$ be a complex manifold,
$T=TX$ its tangent bundle. We have a sheaf of operads ${\cal E}_T$ and a sheaf of PROPs
${\rm END}_T$ on $X$ defined by 
$${\cal E}_T(n) = {\rm Hom}(T^{\otimes n}, T), \quad {\rm END}_T(n,m) = {\rm Hom}
(T^{\otimes n}, T^{\otimes m}).$$
By applying the function $H^\bullet(X, -)$ from sheaves to graded vector spaces,
we get a graded operad
$H^\bullet(X, {\cal E}_T)$ an a graded PROP $H^\bullet(X, {\rm END}_T)$.
Recall also (1.4) that we have the ``big Chern classes"
$\hat c_m(T) \in H^m(X, {\rm Cyc}^m(\Omega^1))$ of the tangent bundle.  Now, a more
inclusive formulation of the properties of the Atiyah class is by using the suspension
of the above PROP and goes
as follows.

\proclaim (3.5.1) Theorem. The element $\Sigma^{-1}\alpha_{TX}\in \Sigma^{-1}H^\bullet
(X, {\cal E}_T)(2)$ is a Lie element. Furthermore, 
the correspondence
$$[x_1,x_2]\in {\rm LIE}(2,1) \quad\mapsto\quad \Sigma^{-1}\alpha_{TX}\in\Sigma^{-1}
H^\bullet(X, {\rm END}_T)(2,1),$$
$$\kappa_n\in {\rm LIE}(n,0) \quad\mapsto\quad \Sigma^{-1}\hat c_m(T)\in \Sigma^{-1}
H^\bullet(X, {\rm END}_T)(n,0),$$
defines a morphism of PROPs
$$ {\rm LIE}\to \Sigma^{-1}
H^\bullet(X, {\rm END}_T) = H^\bullet(X, {\rm END}_{T[-1]}).$$

The proof follows readily from the cohomological Bianchi identity (1.2).

\vskip .3cm

\noindent {\bf (3.6) Weak Lie operad and PROP.}
We denote by ${\cal WL}ie$ the dg-operad governing weak Lie algebras (2.4). It is generated by
elements $\beta_n\in {\cal WL}ie(n)$, ${\rm deg}(\beta_n)=2-n$, $n\geq 2$,
which are antisymmetric with respect to $S_n$ and satisfy the conditions
obtained from Definition 2.4.1.  Thus, the cohomology operad $H^\bullet_d({\cal WL}ie)$
is just ${\cal L}ie$.

\vskip .2cm

The operad ${\cal WL}ie$ can be also described as the cobar-construction of the
commutative operad [GiK]. Explitly, this means that a basis in ${\cal WL}ie (n)$ is formed by
certain trees. More precisely, by an $n$-tree we mean a connected oriented graph
$\Gamma$ with no loops, equipped with structures satisfying the conditions listed below.

\vskip .1cm

\item{(1)} Each vertex of $\Gamma$ has  valency at least 3. In addition, $\Gamma$ has $n+1$ 
legs, i.e., edges bounded by a vertex from one side only.

\vskip .1cm

\item{(2)}  For every vertex $v$ all edges incident to $v$, except exactly one, are oriented
towards $v$. The set of such edges is denoted by ${\rm In}(v)$.

\vskip .1cm

\item{(3)} It follows that all the legs of $\Gamma$ expect exactly one, are oriented towards
$\Gamma$. The set of such legs is denoted by ${\rm In}(\Gamma)$. 

\vskip .1cm

\item{(4)} The set ${\rm In}(\Gamma)$ is identified with $\{1, 2, ..., n\}$. 

\vskip .2cm

Let ${\cal T}(n)$ be the set of isomorphism classes of $n$-trees. For $\Gamma\in {\cal T}(n)$
set
$$\det(\Gamma) = \bigotimes_{v\in {\rm Vert}(\Gamma)} \bigwedge^{\rm max} ({\bf C}^{{\rm In}(v)}).
\leqno (3.6.1)$$

\proclaim (3.6.2) Proposition. We have an identification of 
graded vector spaces,
$${\cal WL}ie(n) = \bigoplus_{\Gamma\in {\cal T}(n)} \det(\Gamma)^*, \quad 
{\rm deg}(\det(\Gamma)^*) = \sum_{v\in {\rm Vert}(\Gamma)} (2-|{\rm In}(v)|),
$$
with the differential being dual to the map given by contraction of edges, 
and the operad structure given by the grafting of trees, see [GiK].

\noindent {\sl Proof:} The identification is obtained by associating to $\beta_n$ the unique $n$-tree
with one vertex (``corolla") and to any composition of the $\beta_n$ the  tree
describing the composition.  The terms in the generalized Jacobi identity correspond, in
geometric language, to all possible $n$-trees with exactly two vertices and one edge
(so that the corolla is obtained from such a tree by contracting this unique edge).

\vskip .2cm

Let WLIE be the (dg-) PROP corresponding to the dg-operad ${\cal WL}ie$ as described in (3.3). It also
has a natural graphical description.  Namely, call an $(n,m)$-graph a (not necessarily connected)
oriented raph $\Gamma$  with $n+m$ legs, of which $n$ are inputs and are labelled with
numbers $1,...,n$, and $m$ are outputs and are labelled by $1, ..., m$, and which satisfy
the conditions (1)-(2) above. Each component of an $(n,m)$ graph is either a tree
satisfying (1)-(3), or a graph with no output. Let ${\cal G}(n,m)$ be the set of
isomorphism classes of $(n,m)$-graphs. Retaining the same notations Vert, In, det,
as for trees, we easily conclude the following.

\proclaim (3.6.3) Proposition.  We have identifications
$${\rm WLIE}(n,m) = \bigoplus_{\Gamma\in {\cal G}(n,m)} \det(\Gamma)^*, \quad 
{\rm deg}(\det(\Gamma)^*) = \sum_{v\in {\rm Vert}(\Gamma)} (2-|{\rm In}(v)|),
$$
with the differential being dual to the map given by contraction of edges, composition
maps given by grafting of graphs, and juxtaposition maps given by disjoint union. 

\vskip .2cm

If ${\cal P}$ is any dg-operad, a family of elements $p_n\in {\cal P}(n)$, $n\geq 2$,
is called {\it a weak Lie    family}, if the  correspondence $\beta_n\mapsto p_n$
gives a morphism of dg-operads ${\cal WL}ie\to {\cal P}$. If $(p_n)$ is a weak Lie
family, then the class of $p_2$ in $H^0({\cal P}(2))$ is a Lie element in $H^\bullet({\cal P})$.
In this case the $p_n$ give also a morphism of PROPs
${\rm WLIE}\to \Pi_{\cal P}$. 

\vskip .3cm

\noindent {\bf (3.7) Differential covariants and the weak Lie PROP.}  We now want to restate 
Theorem 2.6 (which describes the unraveling of the Jacobi identity for
the Atiyah class in the framework of K\"ahler geometry) in a more
universal form.  

Notice, first of all, that the structure we really used, was not the K\"ahler metric itself but
only its canonical (1,0)-connection $\nabla$. So let us call a {\it semiflat manifold}
a pair $(X, \nabla)$ where $X$ is a complex manifold and $\nabla$ is a (1,0)-connection
in $TX$ such that $[\nabla, \nabla]=0$. For such a connection we define $R=[\bar\partial, \nabla]$
and all the considerations of (2.6),  (2.8) hold true.  

Fix $d,i,m,n$ and let ${\cal C}_r$ be the sheaf of semiflat (0,1)-connections on ${\bf C}^r$.
Following Gilkey [Gil] and Epstein [E], introduce the space $V_r^i(n,m)$ of  (not necessarily linear)
differential operators of finite order ${\cal C}\to \Omega^{0,i}\otimes
 {\rm Hom}(T^{\otimes n}, T^{\otimes m})$ defined in some neighborhood of 0, and let 
$U^i_r(n,m)\i V^i_r(n,m)$ be the subspace of operators equivariant under the group
of holomorphic diffeomorphisms. Elements of the latter space will be called differential
covariants of type $(i, n,m)$ of $r$-dimensional semiflat manifolds, since for each such manifold
$(X, \nabla)$ they produce natural tensors in $\Omega^{0,i}({\rm Hom}(T^{\otimes n},
T^{\otimes m}))$. In particular, they do so for each K\"ahler manifold.  In fact, we can say,a
that elements of $U^i_r(n,m)$ are differential covariants of K\"ahler manifolds
which depend only on  the canonical connection. 
The differential
$\bar\partial$ makes $U_r^\bullet(n,m)$ into a complex; taken for all $n,m$, these
complexes form a dg-PROP $U_r^\bullet$.  

For example, $R_n = \nabla^{n-2}R$ 
(the covariant derivative of the curvature) is an
element of $U_r^1(n,1)$. Furthermore, let $\Gamma$ be a $(n,m)$-graph with $N$ vertices.
For every vector space $W$ we have the contraction map
$$p_\Gamma: \bigotimes_{v\in {\rm Vert}(\Gamma)} {\rm Hom}(S^{|{\rm In}(v)|}W, W) 
\to {\rm Hom}(W^{\otimes n}, W^{\otimes m}).
\leqno (3.7.1)$$
Applying this to the tensor product of the $R_{|{\rm In}(v)|}
\in \Omega^{0,1} ({\rm Hom}(S^{|{\rm In}(v)|}T,T))$, we get a covariant
$$R_\Gamma = p_\Gamma\biggl(\bigotimes R_{|{\rm In}(v)|}\biggr)
 \in U^N_r(n,m). \leqno (3.7.2)$$
Because of the symmetry of the $R_i$, the desuspended element $\Sigma^{-1}R_\Gamma$
can be viewed as a morphism
$$\Sigma^{-1}R_\Gamma: \det(\Gamma)^*\to \Sigma^{-1}U^\bullet_r(n,m).\leqno (3.7.3)$$

\proclaim (3.7.4) Theorem. (a) The maps $\Sigma^{-1}R_\Gamma$ define a morphism 
of dg-PROPs $\rho: {\rm WLIE}\to \Sigma^{-1}U_r^\bullet$. \hfill\break
(b) For any $i,n,m$ the morphism of vector spaces
$\rho^i_{n,m}: {\rm WLIE}^i(n,m)\to (\Sigma^{-1}U_r)^i(n,m)$ is surjective.\hfill\break
(c) If $r\gg i,m,n$, then the morphism $\rho^i_{n,m}$ is in fact bijective.

Part (c) means that the ``stabilized" PROP of  differential covariants is just the suspension of the
weak Lie PROP. 

\vskip .2cm

\noindent {\sl Proof:} (a) Follows from Theorem 2.6.

\vskip .1cm

(b) Covariants of $\nabla$ can be viewed as covariants of the total connection $\tilde\nabla = 
\nabla  + \bar\partial$.  It is known classically that all covariants of an affine 
connection are obtained from the covariant derivatives of the curvature (and torsion)
by performing ``tensorial contractions". For example, the argument sketched in [E]
exhibits the Taylor expansion of the Christoffel symbols in the normal coordinates in such a form,
and this clearly suffices. In our case, the covariant derivatives of the curvature of 
$\tilde\nabla$ all reduce to the $R_n$, while a way to perform the contractions produces an $(n,m)$-graph.

 \vskip .1cm

(c) This follows from the main theorem of invariant theory which implies that for $\dim(W)\gg n_1, ... ,
n_N, N$, the space of all $GL(W)$-equivariant maps
$$\bigotimes_{i=1}^N {\rm Hom}(S^{n_i}W, W) \to {\rm Hom}(W^{\otimes n},  W^{\otimes m})$$
has as its basis, the maps $p_\Gamma$ for various $(n,m)$-graphs $\Gamma$ with
 $N$ vertices of valencies $n_i$.

\vfill\eject

\centerline {\bf \S 4. The weak Lie operad in formal geometry and Gelfand-Fuks cohomology.}

\vskip 1cm

\noindent {\bf (4.1) The cochains.}  In this section we describe another way of unraveling the
Jacobi identity for the Atiyah class which uses ``formal geometry" 
(analysis in the space of infinite jets, see [B] [FGG] [GGL]) instead of K\"ahler geometry.
This approach has the advantage of being purely holomorphic.  Instead of Dolbeault cochains,
we will use the following lemma to represent necessary cohomology classes.

\proclaim (4.1.1) Lemma. Let $X$ be a complex manifold and $p:A\to X$ be a locally trivial
fibration with fibers isomorphic to ${\bf C}^N$ for some $N$. Then for any coherent sheaf
$\cal F$ on $X$ we have a natural morphism
$$\tau: \Gamma(A, \Omega^\bullet_{A/X}\otimes p^*{\cal F})\to R\Gamma(X, {\cal F}).$$
If $A$ is a Stein manifold, then $\tau$ is a quasi-isomorphism.

The first statement means that any closed  relative $i$-form on $A$ with values in $p^*{\cal F}$
gives rise to a class in $H^i(X, {\cal F})$. The second statement means that if $A$ is Stein,
then this correspondence is 1-to-1. 

\noindent {\sl Proof:} Follows from the quasiisomorphism
${\cal O}_X\to p_*\Omega^\bullet_{A/X}$ (i.e., from the acyclicity of the global holomorphic
de Rham complex of ${\bf C}^N$).

\vskip .2cm

It was proved by Jouanolou [J] that if $X$ is a quasi-projective algebraic manifold, then
there always exists an $A$ as above which is an affine  variety (therefore Stein). We will
be interested in some natural ${\bf C}^N$-fibrations which, though not Stein in general,
still give the holomorphic cohomology classes we need.

\vskip .3cm

\noindent {\bf (4.2) Formal exponential maps.}
Let $X$ be a complex manifold. Consider the space $\Phi^{(n)}(X)\buildrel p_n\over\rightarrow X$
 of ``$n$th order exponential maps", cf. [B] . By definition, for
$x\in X$ the fiber $\Phi^{(n)}_x(X)$ is the space of $n$th order jets of holomorphic
maps $\phi: T_xX\to X$ such that $\phi(0)=x$, $d_0\phi = {\rm Id}$. Thus $\Phi^{(2)}(X)
=\Phi(X)$ is the affine fibration (2.2) defining torsion-free connections. 
Thus we have a chain of projections
$$ X\leftarrow \Phi^{(2)}(X)\leftarrow \Phi^{(3)}(X) \leftarrow ... \leqno (4.2.1)$$
Each $\Phi^{(n+1)}$ is an affine bundle over $\Phi^{(n)}$ whose associated vector bundle is
\hfill\break
$p_n^* {\rm Hom}(S^{n+1} TX, TX)$. Thus each fiber of $\Phi^{(n)}(X)$
is isomorphic to ${\bf C}^N$ for some $N$ and Lemma 4.1.1 is applicable: every closed relative
form on $\Phi^{(n)}(X)$ gives a holomorphic cohomology class on $X$. 

\vskip .2cm

\noindent {\bf (4.2.2) Example.} Since the space $\Phi(X) = \Phi^{(2)}(X)$
is an affine bundle over \hfill\break
${\rm Hom}(S^2TX, TX)$, it carries a tautological 1-form
$\alpha_2\in \Omega^1_{\Phi(X)/X}\otimes p^* {\rm Hom}(S^2TX, TX)$. This form is relatively closed and
represents, via Lemma 4.1.1,  the Atiyah class $\alpha_{TX}$.

\vskip .2cm

Let $J^{(n)}(TX)\to X$ be the group bundle whose fiber over $x\in X$ is the group
of $n$th jets of  biholomorphisms $\psi: T_xX\to T_xX$ with $\psi(0)=0$,
$d_0\psi = {\rm Id}$. Then $\Phi^{(n)}(X)$ is a bundle of $J^{(n)}(TX)$-torsors. 
Let ${\bf j}^{(n)}(TX)$ be the bundle of Lie algebras associated to $J^{(n)}(TX)$.
Note that we have a natural splitting
$${\bf j}^{(n)}(TX) =\bigoplus_{i=2}^n {\rm Hom}(S^iTX, TX),\leqno (4.2.3)$$
induced by the action of $GL(TX)$ on  ${\bf j}^{(n)}(TX)$.  

Let now $\Phi:=\Phi^{(\infty)}(X) \buildrel p\over\to X$ be the inverse limit of the diagram 
(4.2.1), i.e., the space
of {\it formal exponential maps}. It is a bundle of $J^{(\infty)}(X)$-torsors, where
$J^{(\infty)}(TX) = \lim J^{(n)}(TX)$.  The Lie algebra bundle of the bundle
of proalgebraic groups $J^{(\infty)}(TX)$ is just
$${\bf j}^{(\infty)}(TX) = \prod_{n\geq 2} {\rm Hom}(S^nT, T), \quad T=TX. \leqno (4.2.4)$$
 Denote by $p^{(n)}: \Phi\to \Phi^{(n)}(X)$ the projection and 
set 
$$\Omega^\bullet_{\Phi/X} = \bigcup_n (p^{(n)})^* \Omega^\bullet_{\Phi^{(n)}(X)/X}.
\leqno (4.2.5)$$
As with any bundle of torsors, we have the tautological relative 1-form
$$\omega \in \Omega^1_{\Phi/X}\otimes p^* {\bf j}^{(\infty)}(TX).
\leqno (4.2.6)$$
Projecting $\omega$ to the $n$th graded component in (4.2.4), we get the
{\it tautological form}
$$\alpha_n \in \Omega^1_{\Phi/X} \otimes p^* {\rm Hom}(S^nT, T).\leqno (4.2.7)$$
These forms are formal geometry analogs of the covariant derivatives of the curvature
in (2.5) and satisfy very similar identities, as we shall explain later. 

\vskip .1cm

By the above, for every coherent sheaf $\cal F$ on $X$ 
set
$$A^\bullet_\infty({\cal F}) = a \bigl(\Phi, \Omega^\bullet_{\Phi/X}\otimes p^*{\cal F}
\bigr).$$
This is a complex naturally mapping into $R\Gamma(X, {\cal F})$. Accordingly, for any
coherent sheaf of operads $\cal P$ on $X$ (i.e., an operad in the category of coherent
sheaves) we have a dg-operad $A^\bullet_\infty({\cal P})$. Similarly for PROPs. 

\vskip .1cm

Let us conisder, in particular, the sheaf of operads ${\cal E}_T = \{{\rm Hom}(T^{\otimes n}, T)\}$
and the sheaf of PROPs  ${\rm END}_T = \{{\rm Hom}(T^{\otimes n}, T^{\otimes m})\}$.
The tautological form $\alpha_n$, $n\geq 2$, gives an element
of $A^\bullet_\infty {\cal E}_T(n) \i A^\bullet_\infty {\rm END}_T(n,1)$,
which is antisymmetric and has degree 1. Consider the desuspended dg-PROP
$\Sigma^{-1}A^\bullet_\infty {\rm END}_T$. The shifted tautological forms $\Sigma^{-1}\alpha_n$
becomes antisymmetric of degree $2-n$. Further, let $\Gamma$ be an $(n,m)$-graph
(3.6) with $l$ vertices. We denote by $\alpha_\Gamma \in A^l_\infty {\rm END}_T(n,m)$
the composition of the tautological forms $\alpha_{|{\rm In}(v)|}$ for all vertices $v$ of $\Gamma$,
by using the contractions along the edges of $\Gamma$. Then, because of the antisymmetry
of the $\Sigma^{-1}\alpha_n$, we have that
$$\Sigma^{-1}\alpha_\Gamma \in {\rm Hom}(\det(\Gamma)^*, \Sigma^{-1}A^\bullet_\infty
{\rm END}_T(n,m)).$$ 
Now, a formal geometry version of Theorem 2.6 is as follows. 

\vskip .1cm

\proclaim (4.3) Theorem. The elements $\Sigma^{-1}\alpha_n\in \Sigma^{-1}A^\bullet_\infty
{\cal E}_T(n)$ form a weak Lie family. Moreover, the maps
$$\Sigma^{-1}\alpha_\Gamma:  \det(\Gamma)^* \to 
 \Sigma^{-1}A^\bullet_\infty
{\rm END}_T(n,m), \quad \Gamma \in {\cal G}(n,m),$$
define a morphism of PROPs  ${\rm WLIE}\to \Sigma^{-1}A^\bullet_\infty {\rm END}_T$. 
In particular, the complex of sheaves $p_*\Omega^{\bullet -1}_{\Phi/X}\otimes T$ on $X$
(quasi-isomorphic to $T[-1]$) has a natural structure of a weak Lie algebra. 

To formulate the companion theorem for a vector bundle $E$, we can proceed in a similar way,
by working on the product $\Phi\times_X C$, where $C$ is the fiber bundle over $X$ whose
fiber at $x$ consists of infinite jets of fiberwise linear isomorphisms
$E_x\times X\to E$ identical over $x$. We leave this to the reader.

\vskip .2cm

We know three proofs of Theorem 4.3. The first two will be sketched, and the third one given
in more detail.

\vskip .2cm

\noindent {\bf (4.3.1) First proof (sketch):} We can mimic all features of K\"ahler geometry
but on the space $\Phi$. 
First of all, the bundle $\Phi\to X$ (like other infinite jet bundles, see [FGG] [GGL]), carries a natural
(non-linear) formally integrable connection $D$.  Its covariantly constant sections
over  a simply connected $U\i X$ correspond to affine structures on $U$, i.e., embeddings
of $U$ into an affine space of the same dimension, modulo affine equivalence. 
This decomposes the tangent space $T_\phi\Phi$ at every point
into a direct sum $T_\phi^{1,0}\Phi + T_\phi^{0,1}\Phi$, where
$T_\phi^{0,1}$ is the tangent space to the fiber of $p$ passing through $\phi$ and
$T_\phi^{1,0}$ is the horizontal subspace of the connection.  Accordingly, we have decompositions
$$\Omega^m_\Phi \simeq \bigoplus_{a+b=m} \Omega_\Phi^{a,b}, \quad 
\Omega^{a,b}_{\Phi} = p^*\Omega_X^a\otimes \Omega_{\Phi/X}^b,$$
and the de Rham differential is decomposed as
$$ d=d'+d'', \quad d''=d_{\Phi/X}, \quad (d')^2 = (d'')^2 = [d', d'']=0. $$
We can speak therefore about (0,1) and (1,0)-connections in fiber bundles on $\Phi$.
Every bundle of the form $p^*E$,  lifted from $X$, has a canonical integrable (0,1)-connection. 
The bundle $p^*T$ has, in addition, a natural integrable (1,0)-connection $\nabla$
satisfying the identities
$$[d'', \nabla ] = \alpha_2, \quad \nabla\alpha_n = \alpha_{n+1},$$
which imply our theorem in the same way as in the K\"ahler case.

\vskip .2cm

\noindent {\bf (4.3.2) Second proof (sketch):} 
In line with (2.8), we conisder the odd derivation $D$ of the algebra
$$\Omega^\bullet_{\Phi/X} \otimes p^* \hat S(T^*)$$
obtained by extending $d_{\Phi/X} +\sum_{n\geq 2}\alpha_n^*$. Then we have only to 
prove that $D^2=0$. To do this, we consider, as in (2.8), the fiber bundles
$$\pi: X^{(\infty)}_T\to X, \quad\rho: X^{(\infty)}_{X\times X}\to X$$
where $X^{(\infty)}_T$ is the formal neighborhood of the zero section of $TX$ and 
 $X^{(\infty)}_{X\times X}$ is the formal neighborhood of the diagonal in $X\times X$. The algebra
$\hat S(T^*)$ is just ${\cal O}_{X^{(\infty)}_T}$. The pullback to $\Phi$ of the nonlinear bundle $\rho$ 
possesses an integrable connection along the fibers, which gives rise 
to an algebra differential $\Delta$ in $\Omega^\bullet_{\Phi/X}\otimes p^* {\cal O}_{
X^{(\infty)}_{X\times X}}$ satisfying $\Delta^2=0$. 
On the other hand, on $\Phi$ we have the tautological exponential map which is a nonlinear isomorphism
of fiber bundles
$$ {\rm Exp}:  p^* X^{(\infty)}_T\to p^* X^{(\infty)}_{X\times X},$$
and one can verify that Exp is taking $D$ into $\Delta$, thereby proving the theorem.

\vskip .3cm

\noindent {\bf (4.4) 
 Tautological forms and Gelfand-Fuks cohomology.}
Another way of proof of Theorem 4.3 is to reduce it to  known results 
 about the cohomology of the Lie algebra
of formal vector fields, by making use of the general relationship between this cohomology and
tautological forms.  Let us first recall this relationship [B] [FGG] .

\vskip .1cm

Let $r\geq 1$ be fixed. Denote by $G^{(n)}$ the group of $n$th jets of biholomorphisms
$\phi: {\bf C}^r\to {\bf C}^r$ with $\phi(0)=0$, and by $J^{(n)}\i G^{(n)}$ the normal subgroup
formed by $\phi$ with $d_0\phi={\rm Id}$. So we have an exact sequence
$$1\to J^{(n)}\to G^{(n)}\to GL_r\to 1,\leqno (4.4.1)$$
which, moreover, canonically splits (by considering jets of linear transformations),
 making $G^{(n)}$  a semidirect product.

\vskip .1cm

If $X$ is an $r$-dimensional  complex manifold and $x\in X$,  let $F^{(n)}_x(X)$ be the space of
$n$th jets of biholomorphisms $\phi: {\bf C}^r\to X$ with $\phi(0)=x$. This is a
$G^{(n)}$-torsor.  These torsors unite into a principal $G^{(n)}$-bundle $F^{(n)}(X)
\buildrel q_n\over\to X$
called the bundle of $n$th order frames. The quotient $F^{(n)}(X)/GL_r$
is $\Phi^{(n)}(X)$, the space of $n$th jets of exponential maps from (4.2).

\vskip .1cm

Let ${\bf g} = {\bf h}\oplus {\bf k}$ be a Lie algebra split into a semidirect product of two subalgebras
of which {\bf k} is an ideal. Let $M$ be an {\bf h}-module. Because of the identification
${\bf g}/{\bf k} ={\bf h}$, we can regard $M$ as a {\bf g}-module, and form the relative cochain
complex
$$C^\bullet({\bf g}, {\bf h}, M) = {\rm Hom}(\Lambda^\bullet({\bf g}/{\bf h}), M)^{\bf h}.\leqno (4.4.2)$$
Recall the following standard fact about this complex.

\proclaim (4.4.3) Proposition. If {\bf g}, {\bf h}, {\bf k} are the Lie algebras of connected
Lie groups $G,H,K$ so that $G$ is a semidirect product $HK$, then for every $G$-torsor $P$
we have a natural identification
$$C^\bullet({\bf g}, {\bf h}, M) = \Gamma(P/H, \Omega^\bullet_{P/H}\otimes M)^K.$$

Let us apply this to $G=G^{(n)}, H=GL_r, K=J^{(n)}$. Let ${\bf g}^{(n)}$ be the Lie algebra
of $G^{(n)}$. A representation $M$ of $GL_d$ gives rise, in a standard way, to the
functor from the category of $r$-dimensional vector spaces and their isomorphisms
to the category of vector spaces called the Schur functor and denoted by
$W\mapsto S^M(W)$. In particular, the vector bundle $S^M(TX)$ over $X$ is defined.
We denote by $V$ the standard $r$-dimensional representation of $GL_r$.
Take also $P$ to be (fibers of) the principal $G^{(n)}$-bundle $F^{(n)}(X)\to X$.
We obtain the following.

\proclaim (4.4.4) Proposition. We have a natural identification of complexes of sheaves on $X$
$${\cal O}_X \otimes C^\bullet({\bf g}^{(n)}, {\bf gl}_r, M) \simeq
p_{n*}\biggl(\Omega^\bullet_{\Phi^{(n)}(X)/X}\otimes q_n^*(S^M(TX))\biggr)^{J^{(n)}(TX)}.$$

Note that the tautological forms from (4.2.7) give  global sections of the complex in (4.4.4).
They correspond to $M={\rm Hom}(S^n(V), V)$.

By passing to the limit $n\to\infty$, we consider 
$${\rm Vect}_r^0 = \lim_\leftarrow {\bf g}^{(n)}  = 
\prod_{n\geq 1} {\rm Hom}(S^nV, V),\leqno (4.4.5)$$
the Lie algebra of formal vector fields on ${\bf C}^r$ vanishing at 0. This is a topological Lie algebra
and we will consider its continuous cohomology. 

\vskip .3cm

\noindent {\bf (4.5)  The Lie operad in Gelfand-Fuks cohomology. Results of Fuks.}
Taking, for every $n, m\geq 0$, the relative cochain complex
$$C^\bullet({\rm Vect}^0_r, {\bf gl}_r, {\rm Hom}(V^{\otimes n}, V^{\otimes m})) = 
\Pi^\bullet_r(n,m),\leqno (4.5.1)$$
we get a dg-PROP $\Pi_r^\bullet$.
.Let ${\cal H}_r^\bullet$ be the graded PROP formed
by the cohomology of ${\Pi}_r^\bullet$.  By the above, an element of ${\cal H}_r^i(n)$
gives, for each $r$-dimensional complex manifold $X$, a class in $H^i(X, {\rm Hom}(T^{\otimes n},
T^{\otimes m}))$.  Let
$$a_n \in C^1({\rm Vect}_r^0, {\bf gl}_r, {\rm Hom}(V^{\otimes n}, V)), \quad n\geq 2,$$
be the tautological cochain which associates to a formal vector field its degree $n$ homogeneous
component (lying in ${\rm Hom}(S^nV, V)$).  For any $(n,m)$-graph $\Gamma$ with $N$ vertices
let 
$$a_\Gamma \in C^N({\rm Vect}_r^0, {\bf gl}_r, {\rm Hom}(V^{\otimes n}, V^{\otimes m}))$$
be the cochain obtained by contracting the cup product of the $a_{|{\rm In}(v)|}$,
$v\in {\rm Vert}(\Gamma)$ along the edges of $\Gamma$, cf. (4.2). 
Now, Theorem 4.3 can be reformulated as follows.

\proclaim (4.5.2) Theorem. Let
The maps 
$$\Sigma^{-1} a_\Gamma: \det(\Gamma)^*\to \Pi^N_r(n,m)$$
define a morphism of PROPs ${\rm WLIE}\to \Sigma^{-1}\Pi^\bullet_r$.

In this formulation the theorem follows at once from results of D.B. Fuks [Fuk] who studied
stable cohomology of $\Pi_r^\bullet(m,n)$  (when $r$ is big compared to $m,n$ and the number
of the cohomology) and  identified it with the cohomology of a certain graph
complex.  Translated into our language, his results immediately imply the relation with 
LIE and WLIE. More precisely, we deduce the following fact.

\proclaim (4.5.3) Theorem.  As $r\to\infty$,  each term of the complex $\Pi_r^\bullet(m,n)$
stabilizes, so that we have a limit complex $\Pi^\bullet(m,n)$. Taken for all $m,n$, these
complexes form a dg-PROP, which is isomorphic to $\Sigma^{-1}({\rm WLIE})$.

\noindent {\sl Proof:} The existence of the stabilization and its identification with a
graph complex is completely explicit in [Fuk].  Namely, the space $\Pi_r^N(n,m)$ consists
of  $GL_r$-invariant
antisymmetric continuous maps
$$\bigwedge^N\biggl( \prod_{i\geq 2} {\rm Hom}(S^iV,V)\biggr) \to
 {\rm Hom}(V^{\otimes n}, V^{\otimes m}).
\leqno (4.5.4)$$
Thus
$$\Pi_r^N(n,m) = \bigoplus_{(N_i\in {\bf Z}_+)_{i\geq 2}\atop \sum N_i=N}
{\rm Hom}\biggl( \bigotimes_{i\geq 2} 
\bigwedge^{N_i}{\rm Hom}(S^iV,V), \, {\rm Hom}(V^{\otimes n}, V^{\otimes m})\biggr)^{GL_r}.\leqno 
(4.5.5)$$
 Let $\Gamma$ be an $(n,m)$- graph (3.1) with $N$ vertices. Then we have a natural contraction
map
$$p_\Gamma: \bigotimes_{v\in {\rm Vert}(\Gamma)} {\rm Hom}(S^{|{\rm In}(v)|} (V), V) 
\to {\rm Hom}(V^{\otimes n}, V^{\otimes m}), \leqno (4.5.6)$$
which is obviously invariant.
Moreover, when $r\gg 0$, then by the  main theorem of invariant theory such contraction maps
for various $\Gamma$ provide a basis in the space of all invariant maps. This implies the
stabilization of the $\Pi_r^\bullet(n,m)$. 
Let $N_i(\Gamma)$ be the number of $v\in {\rm Vert}(\Gamma)$ with $|{\rm In}(v)|=i$,
and let
$$t(\Gamma): \bigotimes_{i\geq 2} \bigwedge^{N_i(\Gamma)} {\rm Hom}(S^iV, V) \to
{\rm Hom}(V^{\otimes n}, V^{\otimes m})\leqno (4.5.7)$$
be the antisymmetrization of $p_\Gamma$.
Then
$$t: \det(\Gamma)^*\mapsto t(\Gamma)$$
is the desired isomorphism of complexes ${\rm WLIE}(n,m)\to \Pi^\bullet(n,m)$ of degree $m-n$.
To finish the proof, it remains to identify the composition structure in $\Pi^\bullet$ with that
in $\Sigma^{-1}({\rm WLIE})$, which is straightforward.

\vskip .3cm

\noindent {\bf (4.6) Generalization to other operads.} Theorem 4.5.3 can be 
 straightforwardly generalized to any quadratic Koszul operad $\cal Q$ in the sense
of [GiK]. Namely, let ${\rm Vect}_r^0({\cal Q})$ be the Lie algebra of derivations of
$F_{\cal Q}(r)$, the free $\cal Q$-algebra on $r$ generators. 

\proclaim (4.6.1) Theorem.  Set
$$A_{r,{\cal Q}}^\bullet(m,n) =C^\bullet({\rm Vect}_{r}^0({\cal Q}), {\bf gl}_r, {\rm Hom}(V^{\otimes n}, V^{\otimes m}).$$
Then, as $r\to\infty$, each term of $A_{r,{\cal Q}}^\bullet(m,n)$ stabilizes, and the stable
complexes $A_{{\cal Q}}^\bullet(m,n)$ form a dg-PROP $A_{\cal Q}^\bullet$.
This PROP is isomorphic to $\Pi_{{\bf D}({\cal Q})}$,
 the PROP associated (2.10) to the dg-operad ${\bf D}({\cal Q})$,
the cobar-construction of $\cal Q$. In particular, the graded PROP formed by the cohomology of $A^\bullet_{\cal Q}$
is isomorphic to $\Sigma^{-1} \Pi_{{\cal Q}^!}$ where ${\cal Q}^!$ is the Koszul dual
quadratic operad. 

This statement provides a non-symplectic analog of the result of M. Kontsevich [K1]
describing the stable cohomology of the algebra of hamiltonial vector fields. 
Theorem 4.5.3 corresponds to the case when  ${\cal Q}={\cal C}om$, the operad
describing commutative algebras.

\vskip .3cm

\noindent {\bf (4.7) Example: noncommutativization.}  As we could see before, all the properties
of the Atiyah class, including the detailed unraveling of the Jacobi identity, can be deduced
from the careful study of the non-linear fiber bundle on $X$ whose fiber over $x$ is the
 formal neighborhood of $x$, i.e., the spectrum of the completed local algebra
$\hat{\cal O}_{X,x}$. This algebra is free, i.e., isomorphic to ${\bf C}[[t_1, ..., t_d]]$, $d=\dim(X)$,
but there is no canonical identification, the Atiyah class being an obstruction
to choosing such identifications for all $x$ in a holomorphic way.
Taken for all $x\in X$, the algebras $\hat{\cal O}_{X,x}$ arrange themselves into a sheaf
of complete commutative ${\cal O}_X$-algebras ${\cal O}_{X^{(\infty)}_{X\times X}}$ (functions
on the formal neighborhood of the diagonal), which is locally on $X$ isomorphic to
${\cal O}_X[[t_1, ..., t_d]]$. 

For any commutative ring $R$ let  $R\langle\langle t_1, ..., t_d\rangle\rangle$ be the algebra of non-commutative
formal power series in $t_1, ..., t_d$,  with coefficients in $R$, i.e., the completion 
of the free associative algebra on the $x_i$. 
Now let us make the following 
definition. 

\proclaim (4.7.1) Definition.  Let $X$ be a $d$-dimensional complex manifold. 
A noncommutative structure on $X$ is a 
sheaf of complete associative ${\cal O}_X$-algebras {\bf O} on $X$ which locally on $X$
is isomorphic to ${\cal O}_X\langle\langle t_1, ..., t_d\rangle\rangle$, together with an isomorphism
${\bf O}/[{\bf O},  {\bf O}]\to {\cal O}_{X^{(\infty)}_{X\times X}}$.

In other words, such a structure gives,  for every $x\in X$  a ``non-commutative formal
neighborhood" whose ring of functions is ${\bf O}_x$, the fiber of {\bf O} at $x$. 
These rings are noncanonically isomorphic to ${\bf C} \langle\langle t_1, ..., t_d\rangle\rangle$

\vskip .2cm

\noindent {\bf (4.7.2) Example.}
A natural class of examples of manifolds and, more generally,  stacks
with noncommutative structure is provided by the moduli spaces
of vector bundles (as opposed to more general principal $G$-bundles). Namely, if $E$ is a vector
bundle on an algebraic variety $Z$. Suppose that $H^0(X, {\rm Hom}(E,E))={\bf C}$.
Let $\cal M$ be Kuranishi deformation space of $E$, so that
we have a distinguished point $[E]\in {\cal M}$. Then, by the general principles of deformation
theory [GM], the formal neighborhood of $[E]$ in $\cal M$ is the spectrum of
${\bf H}^0_{Lie}(R\Gamma(Z, {\rm End}(E)))$, the zeroth Lie algebra hypercohomology of the dg-Lie
algebra $R\Gamma(Z, {\rm End}(E))$. Here we regard ${\rm End}(E)$ as a sheaf of Lie
algebras with respect to the bracket $[a,b] = ab-ba$, thereby ignoring a richer structure
of an associative algebra. If we do not ignore this structure, we get an associative dg-algebra structure on
$R\Gamma(Z, {\rm End}(E))$. 
Therefore,  associative algebra hypercohomology
$${\bf H}^0_{Ass} (R\Gamma(Z, {\rm End}(E)))$$
(to be precise, here we mean the Hochshild cohomology with {\bf C} coefficients and the
algebra should be modified so as to get rid of the unity),
 will give us an associative algebra whose
quotient by the commutant maps naturally into the Lie algebra cohomology, i.e., into the
completed local ring of $\cal M$ at $[E]$, and under suitable conditions 
(when the bundle is simple and unobstructed) this is an isomorphism.

\vskip .2cm

\noindent {\bf Remark.} As M. Kontsevich communicated to the author upon readind the manuscript,
he also has had the idea equivalent to Definition 4.7.1 and was aware of Example 4.7.2. 

\vskip .2cm

The considerations of this and the earlier sections revolve, as it is clear from contemplating
Theorem 4.5.3,  around
the Koszul dual pair of operads $({\cal C}om, {\cal L}ie)$:
manifolds are described by commutative algebras of functions, while the curvature data
lead to Lie algebras. So they
 can be generalized to manifolds with a noncommutative
structure, if we consider instead the
  dual pair $({\cal A}ss, {\cal A}ss)$, where ${\cal A}ss$ is the (self-dual) operad governing
associative algebras, see [GiK]. Let us summarize briefly this generalization.

\proclaim (4.7.3) Theorem. Let $X$ be a complex manifold with a noncommutative structure,
$T=TX$ is its usual tangent bundle. Then: \hfill\break
(a) The second order obstruction to global trivialization of the noncommutative formal neighborhoods
is a certain class $\alpha_X\in H^1(X, T\otimes T)$ 
(the noncommutative Atiyah class), whose symmetrization is the usual Atiyah class
$\alpha_{TX}$. \hfill\break
(b) The desuspension $\Sigma^{-1}\alpha_X$, regarded as an element of the
operad $\Sigma^{-1}H^\bullet(X, {\cal E}_T)$, is an associative element, i.e., it defines a morphism of operads
${\cal A}ss\to \Sigma^{-1}H^\bullet(X, {\cal E}_T)$.
In particular, for any sheaf $A$ of commutative ${\cal O}_X$-algebras the shifted
cohomology $H^{\bullet-1}(X, T\otimes A)$ has a natural structure of an associative algebra,
given by $\alpha_X$.\hfill\break
(c) The graded Lie algebra structure on  $H^{\bullet-1}(X, T\otimes A)$ defined by the usual
Atiyah class, is obtained from the associative structure in (b) by the standard formula
$[a,b] = ab \pm ba$. In particular, if $A={\cal O}_X$, then the structure of an associative
algebra on $H^{\bullet-1}(X,T)$ is in fact commutative.

\vfill\eject
    
\centerline{\bf \S 5. The symplectic Atiyah class.}    
    
\vskip 1cm    
    
\noindent {\bf (5.1) Symmetry of the Atiyah class.} Let now $X$ be a complex manifold    
equipped with a holomorphic symplectic structure. Let $\omega \in\Gamma(X, \Omega^2)$    
be the symplectic form. We will identify the tangent bundle $T=TX$ with its dual $T^*$    
by means of $\omega$. After this identification, we can view the Atiyah class    
$\alpha_{TX}$ as an element of $H^1(X, S^2(T)\otimes T)$.    
    
\proclaim (5.1.1) Proposition. The element $\alpha_{TX}$ is totally symmetric, i.e.,    
it lies in the summand $H^1(X, S^3(T))$.    
    
\noindent {\sl Proof:}  Let ${\rm Symp}(X)$ be the sheaf of  connections    
in $TX$ preserving the symplectic form $\omega$.  Since for a symplectic vector    
space $V$ the Lie algebra $sp(V)$ of infinitesimal symplectic transformations    
is identified with $S^2(V)$,  the sheaf ${\rm Symp}(X)$ is, by (1.6),  a torsor    
over $\Omega^1\otimes sp(T) \simeq \Omega^1\otimes S^2(T) \simeq T\otimes S^2(T)$. This shows that    
$\alpha_{TX}$ is symmetric with respect to the permutation of the second and third argument.    
Since it is already symmetric in the first two arguments, the assertion follows.    
    
\vskip .2cm    
    
\noindent {\bf (5.1.2) Remarks.} One can right away exhibit an $S^3(T)$-torsor    
from which $\alpha_{TX}$ is obtained by change of scalars. This is the torsor    
${\rm Symp}_{tf}(X)$ of torsion-free symplectic connections. As in (2.2), it    
can be materialized as the sheaf of sections of the fiber bundle $\Psi(X)\to X$    
whose fiber $\Psi_x(X)$ for $x\in X$ is the space of second jets of holomorphic symplectomorphisms    
$\phi: T_xX\to X$ such  $\phi(0)=x, d_o\phi = {\rm Id}$.  Clearly, $\Psi_x(X)$ is    
an affine space over $S^3(T_xX)$, and sections of $\Psi$ are the same as torsion-free    
symplectic connections.     
    
\vskip .3cm    
    
\noindent {\bf (5.2) The IHX relation for the Atiyah class.} Let $V$ be a finite-dimensional    
 symplectic vector space whose symplectic form is denoted by  $\omega$. Then $V^*$ is also a symplectic vector space, with respect to the inverse form $\omega^{-1}$.    
Let $\Gamma$ be a finite 3-valent graph with possibly several legs (non-compact edges bound by    
a vertex from one side only, cf [GeK]). Denote by ${\rm Vert}(\Gamma), {\rm Ed}(\Gamma),    
{\rm Leg}(\Gamma)$ the sets of vertices, (compact) edges and legs of $\Gamma$.     
Let also ${\rm Flag}(\Gamma)$ be the set of all flags consisting of a vertex and an    
incident half-edge (including a leg) of $\Gamma$. For a vertex $v$ let ${\rm Flag}(v)$    
be the 3-element set of flags having $v$ as a vertex.     
We will distinguish between arbitrary automorphisms of $\Gamma$ and strict    
automorphisms (i.e., those fixing each leg).     

  For a finite-dimensional vector space $W$ we will denote by $\det(W)$ the top
exterior power of $W$. If $I$ is a finite set, then $\det({\bf C}^I)$ will be
abbreviated to $\det(I)$. Note that $\det(I)^{\otimes 2}$ is canonically
(i.e., ${\rm Aut}(I)$-equivariantly) isomorphic to {\bf C}. For an edge
$e$ of $\Gamma$ we denote by ${\rm OR}(e)$ the orientation line of $e$, i.e.,
${\rm OR}(e) = \det(\partial e)$ where $\partial e\i {\rm Flag}(\Gamma)$ is the set
formed by the two flags with edge $e$.

 With these notations, note that
we have a natural $Sp(V)$-equivariant projection    
$$p_\Gamma:  (S^3(V))^{\otimes {\rm Vert}(\Gamma)}\to (V^{\otimes {\rm Leg}(\Gamma)})
\otimes \bigotimes_{e\in {\rm Ed}(\Gamma)} {\rm OR}(e),
\leqno (5.2.1)$$    
obtained by applying the form $\omega: V\otimes V\to {\bf C}$ to any edge of $\Gamma$.   
The factors ${\rm OR}(e)$ appear because of the antisymmetry of $\omega$. 
    
\vskip .2cm    
    
For example, there is a unique, up to scalar, $Sp(V)$-equivariant antisymmetric map    
$$p_{IHX}: S^3(V)\otimes S^3(V)\to S^4(V), \quad p_{IHX}(a\otimes b) = \{a,b\},\leqno (5.2.2)$$    
the Poisson bracket of $a$ and $b$ considered as cubic polynomial functions on $V^*$.    
By working out the definition of the Poisson bracket, we find that the composition    
of $p_{IHX}$ with the embedding $S^4(V)\hookrightarrow V^{\otimes 4}$ can be represented    
as the sum of three projections     
$$p_{IHX} = p_{\bf I} + p_{\bf H} + p_{\bf X}: S^3(V)\otimes S^3(V)\to V^{\otimes 4},\leqno (5.2.3)$$    
where ${\bf I}, {\bf H}, {\bf X}$ are the three possible (up to strict isomorphism) trivalent graphs with    
two vertices and the set of legs identified with $\{1,2,3,4\}$.     
    
\vskip .2cm    
    
Now the Bianchi identity (1.2.2) gives, after the symmetrization, that the Atiyah class    
$\alpha_{TX}\in H^1(X, S^3(T))$ satisfies the so-called IHX relation:    
$$ p_{IHX}(\alpha_{TX}\cup\alpha_{TX}) = \{\alpha_{TX}, \alpha_{TX}\} = 0 \quad {\rm in}    
\quad H^2(X, S^4(T)).\leqno (5.2.4)$$    
Of course, this can be understood from the point of view of the Lie operad, as in (3.5-6),    
the graphs ${\bf I}, {\bf H}, {\bf X}$ corresponding to the three terms in the Jacobi identity.     
    
\vskip .3cm    
    
\noindent {\bf (5.3) Rozansky-Witten classes.}  Let now $\Gamma$ be a trivalent graph    
without legs having $l$ vertices.  Then the projection $p_\Gamma$ from (5.2.1) takes values in {\bf C}.    
By applying it to $(\alpha_{TX})^l     
\in H^l(X, (S^3(T))^{\otimes {\rm Vert}(\Gamma)})$ we get elements    
$$c_\Gamma(X)\in H^l(X, {\cal O})\otimes     
\det ({\rm Vert}(\Gamma))\otimes \bigotimes_e {\rm OR}(e).\leqno (5.3.1)$$    
The factor $\det({\rm Vert}(\Gamma)$ appears because of the anticommutativity
of the multiplication in the cohomology, while the origin of the ${\rm OR}(e)$
was explained in (5.2).  The following lemma shows that the sign factor
in (5.3.1) is the same as the one considered by Rozansky-Witten [RW] and
Kontsevich [K1].

\proclaim (5.3.2) Lemma.  For a trivalent graph $\Gamma$ without edge-loops there is a natural    
(i.e., ${\rm Aut}(\Gamma)$-equivariant) identification of 1-dimensional vector spaces   
$$\det({\rm Vert}(\Gamma))\otimes \bigotimes_{e\in {\rm Ed}(\Gamma)} {\rm OR}(e) \simeq
\det({\rm Ed}(\Gamma)) \otimes \det(H_1(\Gamma, {\bf C})) \simeq$$
$$\simeq \bigotimes_{v\in {\rm Vert}(\Gamma)} \det({\rm Flag}(v)).$$

\noindent {\sl Proof:} We start with the first isomorphism. Note that
$$\det({\rm Ed}(\Gamma))\otimes  \bigotimes_{e\in {\rm Ed}(\Gamma)} {\rm OR}(e)
\simeq \det\biggl( \bigoplus_{e\in {\rm Ed}(\Gamma)} {\rm OR}(e)\biggr),$$
while the consideration of the chain complex of $\Gamma$,
$$\bigoplus_{e\in {\rm Ed}(\Gamma)}{\rm OR}(e)\to {\bf C}^{{\rm Vert}(\Gamma)},$$
gives
$$\det\biggl(\bigoplus {\rm OR}(e)\biggr) \simeq \det({\bf C}^{{\rm Vert}(\Gamma)})\otimes \det    
(H_1(\Gamma, {\bf C})).$$
This implies that the tensor product of the left and the right hand sides of the first
proposed isomorphism in (3.3.2), is canonically trivial.  Because $\det(I)^{\otimes 2}\simeq{\bf C}$
for any finite set $I$,  we get the first isomorphism.

\vskip .1cm

To establish the second isomorphism,  consider the projections    
$${\rm Vert}(\Gamma)\buildrel \phi\over\longleftarrow {\rm Flag}(\Gamma)\buildrel    
\psi\over\longrightarrow {\rm Ed}(\Gamma).$$    
The consideration of fibers of $\psi$ gives    
$$\det({{\rm Flag}(\Gamma}) =    
 \det\biggl( {{\rm Ed}(\Gamma)} \oplus    
\bigoplus_e {\rm OR}(e)\biggr) = \det({{\rm Ed}(\Gamma)}) \otimes \det ({{\rm Vert}(\Gamma)}) \otimes \det(H_1(\Gamma, {\bf C})),$$    
and the consideration of fibers of $\phi$    
gives that    
$$\det({{\rm Flag}(\Gamma)}) = \det ({{\rm Vert}(\Gamma)})\otimes    
\bigotimes_{v\in {\rm Vert}(\Gamma)} \det({{\rm Flag}(v)}),$$    
whence the statement.    
    
\vskip .2cm    
    
For a 3-element set $I$ a choice of direction of the real line $\det({\bf R}^I)$    
is the same as a cyclic order on $I$.    
Thus the classes $c_\Gamma(X)$ can be seen as being elements of $H^l(X, {\cal O})$ but     
defined on graphs with cyclic orders    
on each ${\rm Flag}(v)$ and changing the sign under the changing of the cyclic order.    
Further, it follows from (5.2.3) that the $c_\Gamma$ thus understood satisfy the    
IHX relation in the sense of [RW]. So we get the first part of the following statement:    
    
\proclaim (5.4) Theorem. For any holomorphic symplectic manifold $X$ the classes    
$c_\Gamma(X)$ defined before, give rise to invariants of 3-manifolds with values    
in $H^l(X, {\cal O})$. If $X$ is compact and hyper-K\"ahler, then the $c_\Gamma$ coincide with    
the coefficients  defined by Rozansky and Witten.     
    
The second part just follows from the fact that the curvature represents the Atiyah class     
(Proposition 1.3.1).    
    
It is convenient to consider as in [RW], numerical invariants constructed from    
the $c_\Gamma(X)$. Namely, let us put    
$$\bar c_\Gamma(X) = \omega^{l/2} c_\Gamma(X) \in H^l(X, \Omega^l).\leqno (5.4.1)$$    
Here $\omega\in H^0(X, \Omega^2)$ is the symplectic form and $l$ is the    
(necessarily even) number of vertices of the 3-valent graph $\Gamma$.     
Further, if $X$ is compact and $L$ is a line bundle on $X$, we define the number    
$$b_\Gamma(X, L) = \int_X \bar c_\Gamma(X)\cdot c_1(L)^{\dim(X)-l} \quad\in\quad {\bf C}.    
\leqno (5.4.2)$$    
    
\proclaim (5.4.3) Proposition. The classes $\bar c_\Gamma(X)$ and the numbers $b_\Gamma(X,L)$    
remain unchanged if the symplectic form $\omega$ on $X$ is replaced by $\lambda\omega$,    
$\lambda\in {\bf C}^*$.     
    
\noindent{\sl Proof:} The class $\alpha_{TM}\in H^1(X, S^2T^*\otimes T)$    
does not depend on $\omega$ at all. When we write it in the totally    
symmetric form, we in fact apply the isomorphism    
$$\phi_\omega: S^2T^*\otimes T\to S^2T\otimes T$$    
which is homogeneous in $\omega$ of degree $(-2)$. Since every pairing corresponding to    
an edge of $\Gamma$ is homogeneous in $\omega$ of degree $+1$, we find that    
$c_\Gamma(X)$ is homogeneous of degree    
$$- |{\rm Vert}(\Gamma)| + |{\rm Ed}(\Gamma)| = -(1/2)  |{\rm Vert}(\Gamma)| = -l/2,$$    
where we used the fact that $\Gamma$ is 3-valent. Therefore $\bar c_\Gamma(X)$    
is homogeneous of degree 0.

\vskip .3cm    
    
\noindent {\bf (5.5) Calculation of the $c_\Gamma$ via non-symplectic connections.}     
Note that the Atiyah class $\alpha_{TX}$ used to construct the $c_\Gamma$,    
is defined in terms of the tangent bundle $TX$ alone, without any symplectic    
structure.  Provided such structure $\omega$ is given, $\alpha_{TX}$ just {\it happens to be}    
totally symmetric in all three arguments (if we identify $T\simeq T^*$    
by means of $\omega$). This means that a Dolbeault representative    
of $\alpha_{TX}$ in $\Omega^{0,1}\otimes S^3(TX)$ can be found by forcibly symmetrizing    
the (1,1)-part of the curvature of any (1,0)-connection in $TX$. In particular,    
we can take any K\"ahler (not necessarily hyper-K\"ahler) metric, write its    
curvature as a section of $\Omega^{0,1}\otimes (\Omega^{1,0}\otimes \Omega^{1,0}\otimes T)$,    
identify $T$ with $\Omega^{1,0}$ via the symplectic form and then just symmetrize    
with respect to the last 3 arguments. Denoting by $R$ the corresponding    
$(0,1)$-form with values in $S^3(T)$, we find:    
    
\proclaim (5.5.1) Theorem. The class $c_\Gamma$ is represented by the $(0,l)$-form    
$p_\Gamma(R^{{\rm Vert}(\Gamma)})$.

\vskip .2cm

\noindent {\bf (5.6) Reminder on modular operads and graph complexes.}
 We are now going to upgrade the operadic analysis
of the properties of the Atiyah class and the curvature for general complex manifolds to the
symplectic case. For this, we need the concept of a modular operad [GeK2]. let us briefly recall
this concept.

\vskip .2cm

A stable $((g,n))$-graph is a connected graph $\Gamma$ with the following structures and
properties:

\vskip .1cm

\noindent (5.6.1) The set of legs of $\Gamma$ is identified with $\{1, ..., n\}$.

\vskip .1cm

\noindent (5.6.2) A function $g: {\rm Vert}(\Gamma)\to {\bf Z}_+$ is given such that
$2(g(v)-1)+|v| > 0$ for each vertex $v$.

\vskip .1cm

\noindent (5.6.3) $ \dim (H^1(\Gamma, {\bf C})) + \sum _{v} g(v) = g$.

\vskip .2cm

Let $I((g,n))$ be the set of isomorphism classes of stable $((g,n))$-graphs, 
and $\tilde I((g,n))$ be the similar set in which we allow disconnected graphs as well.

\vskip .1cm

In [GeK2], several versions of modular operads were introduced, differing by the
sign conventions (``cocycles") entering the definition. In the present paper we are going to
use only one of them. Namely, for a stable $((g,n))$-graph $\Gamma$ we set
$$\Re (\Gamma) = \bigotimes_{e\in {\rm Ed}(\Gamma)} {\rm OR}(e). \leqno (5.6.4)$$
The spaces $\Re (\Gamma)$ define a cocycle $\Re$ on the category of graphs in the sense of
[GeK2], \S 4. By  a modular operad we will in the sequel always mean a $\Re$-modular operad.
Explicitly, this is an ordinary operad $\cal P$ with the following additional structures:

\vskip .2cm

\noindent (5.6.5) Symmetry between the inputs and the output, i.e., a $S_{n+1}$-action on ${\cal P}(n)$.
To emphasize this symmetry, we write ${\cal P}((n+1))$ for ${\cal P}(n)$.

\vskip .1cm

\noindent (5.6.6)   A decomposition ${\cal P}((n)) = \bigoplus_{2g-2+n>0}  {\cal P}((g,n))$ into $S_n$-invariant subspaces. The $S_n$-module ${\cal P}((g,n)$ defines, in a standard way,
a functor on the category of $n$-element sets and bijections, whose value on a set $J$
will be denote ${\cal P}((g, J))$.

\vskip .1cm

\noindent (5.6.7) Graphical composition maps
$$\Re (\Gamma)\otimes \bigotimes_{v\in {\rm Vert}(\Gamma)} {\cal P}((g(v), {\rm Leg}(v) \, ))
\to {\cal P}((g,n)).$$

\vskip .1cm

These structures are required to satisfy the compatibility properties given in [GeK2], n. 4.2.

\vskip .3cm

\noindent {\bf (5.7) Examples.}  (a) If $V$ is a vector space with a skew-symmetric
inner product $B$, then we have the {\it endomorphism modular operad}
${\cal E}[V]$ with ${\cal E}[V]((g,n)) =V^{\otimes n}$ for all $g$, the $S_n$-action
being the standard one and the compositions defined by contracting with help of $B$.
Accordingly, for a holomorphic symplectic manifold $X$ we have a graded modular
operad $H^\bullet(X, {\cal E}[T]) = \{ H^\bullet(X, T^{\otimes n})\}$.

\vskip .2cm

(b) The suspended Lie operad $\Sigma {\cal L}ie$ is a modular operad with $\Sigma{\cal L}ie((g,n))=0$ for $g>0$ and
$\Sigma {\cal L}ie(n-1)$  (i.e., the space ${\cal L}ie (n-1)$ placed in degree $n-2$)
for $g=0$, and with the $S_n$-action described in [K1], [GeK1-2].
This action arises naturally from the consideration of Lie algebras with an invariant
inner product. In the same way, $\Sigma{\cal WL}ie$ is a modular dg-operad.

\vskip .2cm

(c) The graph complexes of Kontsevich [K1], once  generalized to allow graphs with legs, form a modular
operad. More precisely, let ${\cal G}((g,n)\i I((g,n))$ be the set of $\Gamma$ for which all the
numbers $g(v)$ are 0. Then the condition (5.6.2) just means $|v|\geq 3$ for all $v$,
and (5.6.3) means that the number of loops in $\Gamma$ is $g$. Set
$${\cal F}((g,n))  = \bigoplus _{\Gamma\in {\cal G}((g,n))} \delta(\Gamma), \quad
\delta(\Gamma) = \det({\rm Ed}(\Gamma))^*\otimes \det(H^1(\Gamma, {\bf C}))^*, \,\,
{\rm deg} (\delta(\Gamma)) = |{\rm Vert}(\Gamma). \leqno (5.7.1)$$
Then the ${\cal F}((g,n))$ form a modular dg-operad $\cal F$  with compositions given by
grafting of graphs and the differential dual to the one contracting edges. In fact,
$\cal F$ is a certain twist of $F{\cal C}om$, the Feynman transform of the commutative
operad defined in [GeK2].  We chose the present version to avoid dealing  in this paper
with twists and suspensions of modular operads and the resulting sign issues. 

Note that the tree part of $\cal F$ is
$${\cal F}((0,n)) = \Sigma{\cal WL}ie ((n)) = \Sigma{\cal WL}ie(n-1). \leqno (5.7.2)$$
Further, ${\cal F}((1,n))$ (the part formed by 1-loop graphs) can be expressed in terms of the
weak Lie PROP, namely it is the subcomplex in
in $\Sigma{\rm  WLIE} (n,1)$ formed by connected graphs.  The legless part ${\cal F}((g,0))$
is the  graph complex  defined by Kontsevich in [K1]. 

\vskip .2cm

(d) Let $\tilde{\cal G}((g,n))$ be the subset of $\tilde I ((g,n))$, see (5.6) formed by graphs with
all $g(v)=0$. Let $\tilde {\cal F}((g,n))$ be the space defined similarly to (5.7.1) but by summing over
$\tilde{\cal G}((g,n))$. They form a modular dg-operad $\tilde{\cal F}$.

\vskip .3cm

\noindent {\bf (5.8) Operadic interpretation of the Jacobi identity: the level of
cohomology.} Let $(X,\omega)$ be a holomorphic symplectic manifold. 
Then we have a graded modular operad $H^\bullet(X, {\cal E}[T])$, see the example
(5.7)(a). 
The Atiyah class $\alpha_{TX}$ can be regarded as
 an element of a modular operad
$$\alpha_{TX} \in H^\bullet(X, {\cal E}[T])((0,3)).$$

\proclaim (5.8.1) Theorem. The correspondence
$$\Sigma( [x_1 ,  x_2]) \in\Sigma  {\cal L}ie(2) = \Sigma {\cal L}ie ((0,3)) \quad \mapsto \quad \alpha_{TX}$$
(with $\Sigma$ meaning the suspension)
defines a morphism of modular operads $\Sigma{\cal L}ie\to  H^\bullet(X, {\cal E}[T])$.

\noindent {\sl Proof:} The only new property here, as compared to Theorem 3.5.1, is that we have a morphism
of modular operads, i.e., that it is invariant  with respect to the action of larger symmetric
groups. But this follows from the total symmetry of $\alpha_{TX}$.

\vskip .3cm

\noindent {\bf (5.9) Operadic interpretation: Dolbeault forms.} We now look at the modular dg-operad
$\Omega^{0,\bullet}({\cal E}[T]) = \{\Omega^{0, \bullet}(T^{\otimes n})\}$.
Assume that $X$ is equipped with a
hyper-K\"ahler metric $h$ . Then the canonical (0,1)-connection $\nabla$ of $h$, see (2.5),
preserves the symplectic structure. The covariant derivatives of the curvature, which we denoted
$R_n, n\geq 2$, are also totally symmetric:
$$R_n\in \Omega^{0,1} (S^{n+1}T). \leqno (5.9.1)$$
For a graph $\Gamma\in \tilde {\cal G}((g,n))$ and a symplectic vector space $W$ let
$$P_\Gamma: \bigotimes_{v\in {\rm Vert}(\Gamma)} S^{|v|}W\to W^{\otimes n} \leqno (5.9.2)$$
be the natural contraction map. Let
$$R_\Gamma = p_\Gamma\biggl(\bigotimes _{v\in {\rm Vert}(\Gamma)} R_{|v|-1}\biggr) \in \Omega^{0,n}(T^{\otimes n}), \,\, N=|{\rm Vert}(\Gamma)|.\leqno (5.9.3)$$
Then  $R_\Gamma$ can be regarded as a morphism
$$R_\Gamma: \delta(\Gamma)\to \Omega^{0, \bullet}(T^{\otimes n}). \leqno (5.9.4)$$

\proclaim (5.9.5) Theorem. The morphisms $R_\Gamma$ extend to a morphism of
modular dg-operads
$$\tilde{\cal F}  \to \Omega^{0,\bullet} ({\cal E}[T]).$$
In particular, for connected graphs with no legs the tensors $R_i$ define a morphism
from ${\cal F} ((g,0))$ (Kontsevich's graph complex) into the Dolbeault complex
$\Omega^{0,\bullet}$. 

\noindent {\sl Proof:} This follows from Theorem 2.6, once we take into account the
additional symmetries of the $R_i$.

\vskip .3cm

\noindent {\bf (5.10) Operadic interpretation: formal geometry.}  Similarly to (4.2), let
$p: \Psi\to X$ be the fiber bundle whose fiber at $x$ consists of all formal {\it symplectic}
exponential maps $T_xX\to X$. This bundle carries the tautological forms
$$\bar\alpha_n \in \Omega^1_{\Psi/X}\otimes p^* S^{n+1}T \leqno (5.10.1)$$
from which we construct the forms
$$\bar\alpha_\Gamma \in \Omega^N_{\Psi/X}\otimes p^*T^{\otimes n}, \quad \Gamma\in \widetilde{\cal
G}((g,n)), |{\rm Vert}(\Gamma)|=N, \leqno (5.10.2)$$
similarly to (5.9.3).

\proclaim (5.10.3) Theorem. The forms $\bar\alpha_\Gamma$ give rise to a morphism of
modular dg-operads
$$\tilde{\cal F}\to H^0(X, (p_*\Omega^\bullet_{\Psi/X})\otimes {\cal E}[T]).$$
In particular,  for connected graphs with no legs these form define a map of 
${\cal F} ((g,0))$ into $p_*\Omega^\bullet_{\Psi/X}$.

The proof can be obtained by embedding $e: \Psi\hookrightarrow \Phi$ where $\Phi$ is the space of
all formal exponential maps from (4.2) and noticing that $\bar\alpha_n = e^*\alpha_n$,
where $\alpha_n$ is the tautological form from (4.2.7). Our statement, which amounts to calculating
$d_{\Psi/X}\bar\alpha_n$, follows from Theorem 4.3. 

\vskip .3cm

\noindent {\bf (5.11) Operadic interpretation: Lie algebra cohomology.} The construction of
(5.10) comes close to the original approach of [K2]: even though we do not
use the $\bar\partial$-foliation on $X$ and its universal characteristic classes, the formal
geometry framework can be regarded as a holomorphic replacement of the 
 $\bar\partial$-theory. In particular, the cohomology of the Lie algebra of formal Hamiltonian vector
fields has direct interpretation in both frameworks. 
The role of Fuks' theorem (4.5.3) from the non-symplectic case is played here  by
the result of Kontsevich  [K1].
Let us formulate it  in a more  general form, allowing graphs with legs so that the operadic
formalism is applicable.

\vskip .2cm

Let $r$ be an even integer and $V={\bf C}^r$ be the standard symplectic vector space of
dimension $r$. Let
$${\rm Ham}_r^0 = \prod_{n\geq 2} S^nV \leqno (5.11.1)$$
be the Lie algebra of formal hamiltonian vector fields on $V$ with trivial constant term. 
Its degree 2 part is $S^2V = {\bf sp}_r$, the Lie algebra of linear symplectomorphisms. 
As in (4.4), any relative cochain of $({\rm Ham}_r^0, {\bf sp}_r)$ with coefficients
in some tensor power of $V$ gives rise to a natural relative form on $\Psi/X$
with values in the corresponding tensor power of $p^*TX$. In particular, the tautological
cochain
$$\bar a_n\in C^1({\rm Ham}^0_r, {\bf sp}_r, S^{n+1}V) \leqno (5.11.2)$$
associating to a vector field its degree $n+1$ part, corresponds to the form $\bar\alpha_n$.
For any graph $\Gamma\in {\cal G}((g,n))$
define
$$\bar a_\Gamma\in C^N({\rm Ham}^0_r, {\bf sp}_r, S^{n+1}V) \leqno (5.11.3)$$
as in (5.8.9). Furthermore, the complexes
$$P^\bullet_r((g,n)) :=  C^\bullet({\rm Ham}^0_r, {\bf sp}_r, S^{n}V) \leqno
(5.11.4)$$
define a modular dg-operad $P^\bullet_r$. The symplectic analog of Theorem 4.5.2
 is as follows.

\proclaim (5.11.5) Theorem. (a) The maps
$$\bar a_\gamma: \delta(\Gamma) \to P_r^\bullet((g,n)), \quad \Gamma\in\tilde{\cal G}((g,n))$$
define a morphism of modular operads $\tau: \tilde{\cal F}\to P^\bullet_r$.
\hfill\break
(b) If $r\gg g,n,i$, then the map of vector spaces
$$\tau: \tilde{\cal F}^i((g,n))\to P^i_r((g,n))$$
is an isomorphism.

This is proved in the same way as the result in [K1] (Theorem 1.1, which concerns graphs without legs, i.e.,
cohomology with trivial coefficients) or [Fuk].

\vfill\eject

\centerline {\bf References}    
    
\vskip 1cm    

\noindent [Ad] J.F. Adams, Infinite Loop Spaces, Princeton Univ. Press, 1977.

\vskip .2cm

\noindent [AL] B. Ang\'eniol, M. Lejeune-Jalabert, Le th\'eor\` eme de Riemann-Roch singulier
pour les $\cal D$-modules, {\it Ast\' erisque} {\bf 130} (1985),  130-160. 

\vskip .2cm
    
\noindent [At] M. F. Atiyah, Complex analytic connections in fiber bundles,  {\it Trans. AMS}, {\bf 85}    
(1957), 181-207 (Shorter version reprinted in his Collected Papers, vol.1, p. 95-102,     
Clarendon Press, Oxford, 1988).    

\vskip .2cm

\noindent [BCOV]  M. Bershadsky, S. Cecotti,  H. Ooguri, C. Vafa, Kodaira-Spencer theory
of gravity and exact results for quantum string amplitudes, {\it Comm. Math. Phys.}
{\bf 165} (1994), 311-428.
    
\vskip .2cm    

\noindent [B] R. Bott, Some aspects of invariant theory in differential geometry, in:
``Differential operators on Manifolds" (CIME),  p. 49-145,  Edizione Cremonese, 1975
(Reprinted in his Collected Papers, vol. 3, p.  357-453, Birkhauser, Boston, 1996). 

\vskip .2cm

\noindent [E] D.B.A. Epstein, Natural tensors on Riemannian manifolds, {\it J. Diff. geom.}
{\bf 10}  (1975), 631-645. 

\vskip .2cm

\noindent [Fuk] D.B. Fuks,  Stable cohomologies of a Lie algebra of formal vector fields
with tensor coefficients, {\it Funct. Anal. Appl.} {\bf 17} (1983), 295-301.

\vskip .2cm

\noindent [FGG] D.B. Fuchs, I.M. Gelfand, A.M. Gabrielov,  The Gauss-Bonnet theorem and
the Atiyah-Patodi-Singer functionals for the characteristic classes of foliations, 
{\it Topology}, {\bf 15} (1976), 165-188.
(Reprinted in: Collected Papers of I.M. Gelfand, Vol. 3, p. 379-402, Springer-Verlag,  1989)

\vskip .2cm

\noindent [GGL] A.M. Gabrielov, I.M. Gelfand, M.V. Losik, Combinatorial computation of
characteristic classes, {\it Funct. Anal. Appl.} {\bf 9} (1975), 103-115.
(Reprinted in: Collected Papers of I.M. Gelfand, Vol. 3, p. 407-419, Springer-Verlag,  1989)

\vskip .2cm    
    
\noindent [GeK1] E. Getzler, M. Kapranov, Cyclic operads and cyclic homology, in:    
``Geometry, Topology and Physics for R. Bott" (S.-T. Yau, Ed.) p. 167-201, International    
Press, Boston, 1995.     
    
\vskip .2cm    
    
\noindent [GeK2] E. Getzler, M. Kapranov, Modular operads, preprint dg-ga/9408003,
to appear in {\it Compositio Math.}

\vskip .2cm

\noindent [Gil] P.B. Gilkey, Curvature and eigenvalues of the Dolbeault complex for K\"ahler manifolds,
{\it Adv. in Math.} {\bf 11} (1973), 311-325. 
    
\vskip .2cm

\noindent [GiK] V. Ginzburg, M. Kapranov, Koszul duality for operads, {\it Duke Math. Journal},     
{\bf 76} (1994), 203-272.     

\vskip .2cm

\noindent [GM] W.M. Goldman, J.J. Millson, The deformation theory of representations of fundamental groups 
of compact K\"ahler manifolds, {\it Publ. IHES}, {\bf 67} (1988), 43-96.
    
\vskip .2cm    

\noindent [J]  J.-L. Jouanolou, Une suite exacte de Mayer-Vietoris dans la K-theorie
algebrique, in: Lecture Notes in Math., {\bf 341}, p. 293-316, Springer-verlag, 1973.

\vskip .2cm

\noindent [K1] M. Kontsevich, Formal noncommutative symplectic geometry, in ``Gelfand Mathematical    
Seminars 1990-92" (L. Corwin, I. Gelfand, J. Lepowski Eds.)  p.  173-187, Birkhauser, Boston, 1993.     
    
\vskip .2cm    
    
\noindent [K2] M. Kontsevich, letter to V. Ginzburg, Jan. 8, 1997.     

\vskip .2cm

\noindent [KM] I. Kriz, J.P. May, Operads, Algebras, Modules and Motives, {\it Asterisque} {\bf 
233}, Soc. Math. France, 1995.

\vskip .2cm

\noindent [M] M. Markl, Models for operads, {\it Comm. in Algebra}, {\bf 24}  (1996), 1471-1500.

\vskip .2cm

\noindent  [RW]  L. Rozansky, E. Witten, Hyper-K\"ahler geometry and invariants of 3-manifolds,    
preprint hep-th/9612216.   

\vskip .2cm

\noindent [S] J. Stasheff, Differential graded Lie algebras, quasi-Hopf algebras and higher
homotopy algebras, in: Lecture Notes in Math., {\bf 1510}, p. 120-137,
Springer-Verlag, 1992. 
    
\vskip .2cm    
    
\noindent [W] R.O. Wells, Differential calculus on complex manifolds, Springer-Verlag, 1982.

\vskip 2cm    
    
\noindent {\sl Department of Mathematics, Northwestern University, Evanston IL 60208, \hfill\break    
email: kapranov@math.nwu.edu}

  \bye